\documentclass[journal]{IEEEtran}

\usepackage{amsfonts}
\usepackage{amssymb}
\usepackage{amsmath}
\usepackage{graphicx}
\usepackage{epstopdf}
\usepackage{cite}
\usepackage{color}
\usepackage{caption,subcaption}
\usepackage{overpic}

\def\scr#1{{\cal #1}}
\newcommand{\R}{\mathbb{R}}
\def\eq#1{\begin{equation}#1\end{equation}}
\def\rep#1{(\ref{#1})}
\newcommand{\bbb}{\mathbb}
\newtheorem{theorem}{Theorem}
\newtheorem{lemma}{Lemma}
\newtheorem{remark}{Remark}
\newtheorem{proposition}{Proposition}
\newtheorem{corollary}{Corollary}
\newtheorem{assumption}{Assumption}

\def\qed{ \rule{.08in}{.08in}}

\newcommand{\1}{\mathbf{1}}
\newcommand{\0}{\mathbf{0}}

\newcommand{\D}{\mathcal{D}}

\newtheorem{definition}{Definition}

\begin{document}

\title{Analysis and Control of a Continuous-Time Bi-Virus Model}

\author{Ji Liu,
        Philip E. Par\'{e},
        Angelia Nedi\'{c},
        Choon Yik Tang,
        Carolyn L. Beck,
        Tamer Ba\c{s}ar
\thanks{Some of the material in this paper was presented at the 55th
IEEE Conference on Decision and Control \cite{liu2016onthe}.}
\thanks{This work is based on research partially sponsored by the National Science Foundation grants ECCS 15-09302, CCF 11-11342, DMS 13-12907, and CNS 15-44953, the Office of Naval Research (ONR) MURI Grant N00014-16-1-2710, US Army Research Office (ARO) Grant W911NF-16-1-0485, and ONR Basic Research grant Navy N00014-12-1-0998.  All material in this paper represents the position of the authors and not necessarily that of the funding agencies.}
\thanks{Ji Liu is with Stony Brook University (\texttt{ji.liu@stonybrook.edu}).
Philip E. Par\'{e},
Carolyn L. Beck, and Tamer Ba\c sar are with the Coordinated Science Laboratory at the University of Illinois at Urbana-Champaign (\texttt{\{philpare, beck3, basar1\}@illinois.edu}).
Angelia Nedi\'{c} is with the School of ECEE  at Arizona State University
(\texttt{angelia.nedich@asu.edu}).
Choon Yik Tang is with the School of ECE  at the University of Oklahoma
(\texttt{cytang@ou.edu}). }
}

\maketitle

\begin{abstract}
This paper studies a distributed continuous-time bi-virus model in which
two competing viruses spread over a network consisting of multiple groups of individuals.
Limiting behaviors of the network are characterized by
analyzing the equilibria of the system and their stability.
Specifically, when the two viruses spread over possibly different directed infection graphs,
the system may have (1) a unique equilibrium, the healthy state, which is globally stable,
implying that both viruses will eventually be eradicated,
(2) two equilibria including the healthy state and a dominant virus state, 
which is almost globally stable, implying that one virus will pervade the entire network causing a single-virus epidemic
while the other virus will be eradicated, or 
(3) at least three equilibria including the healthy state
and two dominant virus states, depending on certain conditions on the healing and infection rates.
When the two viruses spread over the same directed infection graph,
the system may have zero or infinitely many coexisting epidemic equilibria,
which represents the pervasion of the two viruses.
Sensitivity properties of some nontrivial equilibria are investigated in the context of a decentralized control technique,
and an impossibility result is given for a certain type of distributed feedback controller.
\end{abstract}


\IEEEpeerreviewmaketitle

\section{Introduction}

The spread of epidemic processes over large populations is an important research topic,
and is in fact a widely studied one in epidemiology \cite{book}. To model such a process,
various epidemic models have been proposed such as the susceptible-infected-recovered (SIR), susceptible-exposed-infected-recovered (SEIR), and susceptible-infected-susceptible (SIS) models \cite{sir,seir,yorke,sis}.
Bernoulli developed one of the first known models inspired by the smallpox virus \cite{bernoulli1760essai}.
The first SIS model was introduced in \cite{kermack1932contributions}.
In this paper, we focus on the study of distributed SIS epidemic models, where there are two ways to consider such a system: 1) the model consists of $n>1$ interacting individuals
and the evolution of the probability of each individual being infected is studied, or 2) the model consists of $n>1$ groups of individuals
and the evolution of the percentage of infected members of each group is studied.
The first type of SIS model has been studied in both discrete-time 
\cite{WangSRDS03,WangTISSEC08,BokharaieMTNS10,AhnCDC13,AhnThesis}
and continuous-time   \cite{OmicTN09,OmicINFOCOM09,Omicarxiv,KhanaferITA14,KhanaferACC14,KhanaferCDC14,KhanaferAutomatica14,pare2015stability,pare2017epidemic}.
The first multi-agent, probability-based, continuous-time model was proposed by
Van Mieghem et al. \cite{OmicTN09} in which the underlying neighbor graph
is assumed to be undirected. The same model on a directed neighbor graph
has been recently studied by Khanafer et al.~\cite{KhanaferAutomatica14}
for both strongly and weakly connected neighbor graphs.
The second, or group-based type of models, has been studied in \cite{FallMMNP07}.

The idea of competing SIS virus models is  pursued in \cite{nowak1991evolution,prakash2012winner,sahneh2014competitive,santos15,santos2015bivirus,wei2013competing,watkins2015optimal}.
This work is motivated by the competition of different viral strains \cite{nowak1991evolution}, where there are two competing viruses in a human contact network.
These models have a wide range of other applications, including social networks, where the goal 
is to understand how competing opinions spread on different social networks \cite{sahneh2014competitive}, competing products in a market, and agents' opinions about politicians from opposing parties \cite{trpevski2010model}.
Competing SIS models were first introduced in \cite{nowak1991evolution}, which is an extension of \cite{kermack1932contributions}, where the model considers the dynamics of three groups: 1) susceptible, 2) infected with virus one, and 3) infected with virus two. These dynamics are modeled by three differential equations where full connectivity of the agents is assumed
(i.e., the infection graph is a complete graph), and it is also assumed that the two viruses are both homogeneous.\footnote{
We say that a virus is \textit{homogeneous} if all agents have the same infection rate and healing rate.
Otherwise, the virus is called \textit{heterogeneous}.
}
In \cite{prakash2012winner}, two competing homogeneous viruses spreading over the same nontrivial (not necessarily fully connected), undirected, connected network is studied. The set of equilibrium points is determined and sufficient conditions for local stability are given for all equilibria except the coexisting equilibrium.
In \cite{sahneh2014competitive}, the equilibria of two competing homogeneous virus models over the same as well as different undirected graph structures, are studied. Existence of the coexisting epidemic states, where both viruses are at nontrivial (nonzero) equilibria, is shown, but no stability analysis is provided.
Note that all this previous work is conducted for homogeneous viruses over undirected graph structures with limited/local stability analysis.
The following are the two exceptions:
in \cite{santos2015bivirus}, a sufficient condition for the global asymptotic survival of a single virus is given for a model of two competing viruses, both homogeneous in the healing rate and propagating over undirected, regular graphs.
In \cite{watkins2015optimal,watkins2016optimal}, a necessary and sufficient condition for local exponential stability of the origin is provided for two competing heterogeneous viruses over strongly connected graphs. In addition, a geometric program is formulated, working toward optimal stabilization and rate control of the virus. However, stability of the epidemic equilibria (i.e., nonzero equilibria) is not explored.
Note that none of the existing work considers  heterogeneous viruses over directed graph structures and performs global stability analysis, exploring all of the system's equilibria.

Competing viruses are also explored for an SIR 
model in \cite{yagan2013conjoining}. Additionally, recently multiple competing viruses, or multi-virus models, have been explored in  \cite{xu2012multi,acc_multi,arxiv_multi}. A centralized control technique for multi-virus systems is explored in \cite{acc_multi,arxiv_multi}.

Various  control strategies have been explored for the single-virus model \cite{KhanaferAutomatica14,PreciadoTCNS14,bullo2014control}.
In \cite{KhanaferAutomatica14}, an antidote control technique is proposed.
In \cite{PreciadoTCNS14}, an optimal vaccination control technique is developed using geometric programming ideas.
In \cite{bullo2014control}, a network control scheme is applied to a discretized, linearized version of the single-virus model.
To the best of our knowledge, the only control strategies that have been designed for the bi-virus model are the ones in \cite{watkins2015optimal} and \cite{watkins2016optimal}, which use geometric programs and are centralized.

In this paper, we study a distributed continuous-time bi-virus model over {\em directed} graphs.
The model describes how two competing SIS viruses spread over a network of $n>1$ groups of individuals.
By \textit{competing} we mean that no individual can be infected with the two viruses simultaneously.
An individual may be infected with one of the two viruses by individuals in its own group, as well as nearest-neighbor groups.
These neighbor relationships among the $n$ groups are described
by a directed graph on $n$ vertices with an arc (or a directed edge) from vertex $j$
to vertex $i$ whenever the individuals in group $i$ can be infected by those in group $j$.
Thus, the neighbor graph has self-arcs at all $n$ vertices, and the arc direction
represents the direction of the 
contagion.
The two viruses may spread over different infection graphs with each infection graph
being a spanning subgraph of the neighbor graph. Thus,
the neighbor graph is the union of the two infection graphs.\footnote{
A \textit{spanning subgraph} is a subgraph that contains all the vertices of the original graph.
The \textit{union} of two directed graphs with the same vertex set is a directed graph
with the same vertex set whose arc set is the union of the arc sets of the two graphs.}

For two competing viruses, the SIS model is probably the simplest one, but it has limiting behaviors that are already complicated and challenging to analyze. As we will show, the bi-virus SIS model can predict much richer spreading phenomena compared with the single-virus SIS model.
The multi-group bi-virus model can be applied to a number of areas, allowing one to understand dynamics of, for example, two competing products in an economic market \cite{armington}, two competing memes in a social network \cite{meme}, and two competing species in an ecological environment \cite{lion}. Thus, there is ample motivation to thoroughly understand all possible limiting behaviors of the model, which admit different interpretations in different fields. For example, in an epidemic network, a globally stable healthy state, dominant virus state, and coexisting epidemic state predict the eradication of both of the viruses, the triumph of one virus over the other, and the pervasion of both viruses, respectively.
Since different applications may prefer different limiting behaviors, it is of great interest to propose efficient control techniques for promoting and/or precluding the spread process in a competitive environment, preferably in a distributed manner. 
Apparently, for a bi-virus epidemic network, pervasion of any virus is undesirable, and thus a natural question is how to efficiently attenuate or eliminate epidemic spreading using a distributed controller.

The main contributions of this paper are three-fold.
First, we analyze the equilibria of the bi-virus model over directed graphs
and their stability for both homogeneous and heterogeneous viruses.
Second, we derive a sensitivity result for the nontrivial equilibria with respect to the infection and the healing rates,
which demonstrates the effects of the simplest local control (i.e., an individual locally adjusts his/her healing and infection rates).
Third, we provide an interesting and surprising impossibility result for a certain type of distributed feedback controller, which reveals why distributed control of (bi-virus) epidemic networks is a challenging problem. 
All the results are validated by a set of illustrative simulations.

The model and assumptions considered in this paper are more general than
those in the existing literature \cite{prakash2012winner,sahneh2014competitive,santos2015bivirus,santos15,wei2013competing,watkins2015optimal,watkins2016optimal}
in three aspects.
First, the work of \cite{prakash2012winner,sahneh2014competitive,santos2015bivirus,santos15,wei2013competing}
only considers undirected graphs, whereas this paper studies directed graphs.
Second, the work of \cite{prakash2012winner,wei2013competing}
only considers homogeneous viruses, whereas this paper studies heterogeneous viruses.
Third, the analyses in \cite{prakash2012winner,sahneh2014competitive,santos2015bivirus,santos15,wei2013competing,watkins2015optimal,watkins2016optimal}
are limited to local stability of the equilibria, while the majority of the analysis
 performed in this paper is global.
Moreover,
this paper explores two possible distributed control techniques for the bi-virus model, an endeavor that has not been widely pursued before. 

As we will see shortly, the bi-virus model includes the continuous-time single-virus SIS model
as a special case. A byproduct of this paper is thus an analysis of the single-virus system
which has been studied earlier in \cite{Gonorrhea,FallMMNP07,KhanaferAutomatica14}.
The analysis herein is performed under weaker  assumptions on the infection rates $\beta_{ij}$
and healing rates $\delta_{i}$, therefore generalizing previous results on the single-virus system.
Although the system defined in \cite{KhanaferAutomatica14} admits the same
mathematical expression as the single-virus model considered herein, a key difference between the work
of \cite{KhanaferAutomatica14} and this paper, other than different physical meaning of the models,
is that it is assumed in \cite{KhanaferAutomatica14} that $\beta_{ii}=0$,
for all $i$, whereas this is not the case here. Another difference is that it is also assumed in \cite{KhanaferAutomatica14}
that $\beta_{ij}=\beta_{kj}$ for all $i,j,k$ if they are both nonzero.
In \cite{Gonorrhea}, it is assumed that if $\beta_{ij}> 0$, then $\beta_{ji}> 0$ for all $i,j$, though they are not necessarily equal.  
Moreover, in \cite{Gonorrhea,KhanaferAutomatica14,FallMMNP07}, it is assumed that $\delta_i>0$ for all $i$.
In contrast, this paper does not impose any of the aforementioned assumptions and,
thus, considers a more general model than those in \cite{FallMMNP07,KhanaferAutomatica14}. Note that the generalization of allowing $\delta_i=0$ includes a susceptible-infected (SI) model.




Some of the material in this paper was presented in preliminary form in \cite{liu2016onthe};
 this paper presents a more comprehensive treatment of the work in \cite{liu2016onthe}.
Additional contributions of this paper, that are not in \cite{liu2016onthe}, include 1)  
 complete proofs of all the results, 
2)  several extensions, including the establishment of 
uniqueness of the parallel equilibria in Theorems \ref{parallel1} and \ref{parallel2}, 3) viewing the sensitivity analysis from a control perspective, 
4) an impossibility result
for a certain type of distributed feedback controller,
and
5) an in-depth set of simulations, illustrating the results and some unproven phenomena. 



The remainder of the paper is organized as follows.
The
basic properties and assumptions of the system model are given in Section~\ref{system}, 
and  the full probabilistic, $3^n$ state model is presented.
The system equilibria and their stability are studied in Section~\ref{equilibria}. 
The sensitivity of the equilibria is investigated in Section~\ref{sense}.
An impossibility result for distributed feedback control is provided in Section~\ref{feedback}.
Simulations are given in Section~\ref{simulations}.
The paper concludes with some remarks in Section~\ref{conclusion}.
The proofs of some assertions in the paper are given in the appendix.
In the rest of this section, we introduce some notation and provide a number of preliminary results. 

\subsection{Notation}

For any positive integer $n$, we use $[n]$ to denote the set $\{1,2,\ldots,n\}$.
We view vectors as column vectors. We use $x'$ to denote the transpose of a vector $x$ and,
similarly, we use $A'$ for the transpose of a matrix $A$.
The $i$th entry of a vector $x$ will be denoted by $x_i$.
The $ij$th entry of a matrix $A$ will be denoted by $a_{ij}$. 
We use $\0$ and $\1$ to denote the vectors whose entries all equal to $0$'s or $1$'s, respectively,
and $I$ to denote the identity matrix,
where the dimensions of the vectors and matrices are to be understood from the context.
For any vector $x\in\R^n$, we use ${\rm diag}(x)$ to denote the $n\times n$ diagonal matrix
whose $i$th diagonal entry equals $x_i$.
For any two sets $\scr{A}$ and $\scr{B}$,
we use $\scr{A}\setminus \scr{B}$ to denote the set of elements that are in $\scr{A}$ but not in $\scr{B}$, and $\scr{A} \subset \scr{B}$ to denote equality or a proper subset.
The notation $1_{a=b}$ is used as an indicator function which takes value one if $a$ equals $b$ and zero otherwise. For $1_{A=b}$, where $A$ is a matrix and $b$ is a scalar, the result is a binary matrix of the same dimensions as $A$ with entries $1_{a_{ij}=b}$.

For any two real vectors $a,b\in\R^n$, we write $a\geq b$ if
$a_{i}\geq b_{i}$ for all $i\in[n]$,
$a>b$ if $a\geq b$ and $a\neq b$, and $a \gg b$ if $a_{i}> b_{i}$ for all $i\in[n]$.
Similarly, for any two real matrices $A,B\in\R^{m\times n}$, we write $A\geq B$ if
$a_{ij}\geq b_{ij}$ for all $i\in[m]$ and $j\in[n]$,
$A>B$ if $A\geq B$ and $A\neq B$, and $A \gg B$ if $a_{ij}> b_{ij}$ for all $i\in[m]$ and $j\in[n]$.

For a complex number $x$, we use $|x|$ and ${\rm Re}(x)$ to denote its
magnitude and real part, respectively.
For a real square matrix $M$, we use $\rho(M)$ to denote its spectral radius
and $s(M)$ to denote the largest real part among its eigenvalues, i.e.,
\begin{eqnarray*}
\rho(M) &=& \max \left\{|\lambda|\ : \ \lambda\in\sigma(M)\right\}, \\
s(M) &=& \max \left\{{\rm Re}(\lambda)\ : \ \lambda\in\sigma(M)\right\},
\end{eqnarray*}
where $\sigma(M)$ denotes the spectrum of $M$.


The sign function of a real number $x$ is defined as follows:
$$
{\rm sgn}(x) = \left\{ \begin{array}{lll}
  -1 &\mbox{ if \ $x<0$}, \\
  0 &\mbox{ if \ $x=0$}, \\
  1 &\mbox{ if \ $x>0$}.
       \end{array} \right.
$$
Note that for any real number $x\neq 0$,
$\frac{d|x|}{dx} = {\rm sgn}(x)$.

\subsection{Preliminaries}

For any two nonnegative vectors $a$ and $b$ in $\R^n$ ($a,b\geq \0$),
we say that $a$ and $b$ have the same sign pattern if they have zero entries
and positive entries in the same places, i.e., for all $i\in[n]$,
$a_i=0$ if and only if $b_i=0$, and $a_i>0$ if and only if $b_i>0$.
A square matrix is called {\em irreducible}  if it cannot be permuted to a block upper triangular matrix.

\begin{lemma}
Suppose that $Mx=y$ where $M\in\R^{n\times n}$ is an irreducible nonnegative matrix
and $x,y>\0$ are two vectors in $\R^n$.
If $x$ has at least one zero entry, then $x$ and $y$ cannot have the same sign pattern.
In particular, there exists an index $i\in[n]$ such that $x_i=0$ and $y_i>0$.
\label{pattern}\end{lemma}

A real square matrix is called {\em Metzler} if its off-diagonal entries are all nonnegative.
Thus, any nonnegative matrix is Metzler.

\begin{lemma}\label{pos}
For any matrix $M$ and any real number $\phi$, if $A:=M-\phi I$, then $\sigma(M)=\sigma(A)+\phi$.
\end{lemma}

The proof of this lemma is simple and therefore omitted.

The following results from Chapter 2 of \cite{varga} for nonnegative matrices, which also hold for Metzler matrices by Lemma \ref{pos}, with $\phi = \min\{0,m_{11},...,m_{nn}\}$, will be used in the subsequent analysis.

\begin{lemma}
{\rm (Lemma 2.3 in \cite{varga})}
Suppose that $M$ is an irreducible Metzler matrix. Then,
$s(M)$ is a simple eigenvalue of $M$ and there exists a unique (up to scalar multiple) vector $x\gg \0$ such that
$Mx=s(M)x$.
\label{metzler0}\end{lemma}

\begin{lemma}
{\rm (Section 2.1 in \cite{varga})}
Suppose that $M$ is an irreducible Metzler matrix in $\R^{n\times n}$ and $x>\0$ is a vector in $\R^n$.
For any $ \lambda \in \R$, if $Mx < \lambda x$, then $s(M)<\lambda$.
If $Mx = \lambda x$, then $s(M)=\lambda$. If $Mx > \lambda x$, then $s(M)>\lambda$.
\label{metzler}\end{lemma}

We now introduce the Perron-Frobenius Theorem for irreducible nonnegative matrices.

\begin{lemma}
{\rm (Theorem 2.7 and Lemma 2.4 in \cite{varga})}
Suppose that $M$ is an irreducible nonnegative matrix. Then,
\begin{enumerate}
    \item $M$ has a positive real eigenvalue equal to its spectral radius, $\rho(M)$.
    \item $\rho(M)$ is a simple eigenvalue of $M$.
    \item There is an eigenvector $x\gg \0$ corresponding to $\rho(M)$.
    \item $\rho(M)$ increases when any entry of $M$ increases.
    \item If $N$ is also an irreducible nonnegative matrix and $M > N$, then $\rho(M) > \rho(N)$.
\end{enumerate}
\label{perron0}\end{lemma}



\begin{proposition}
Suppose that $\Lambda$ is a negative diagonal matrix in $\R^{n\times n}$
and $N$ is an irreducible nonnegative matrix in $\R^{n\times n}$.
Let $M=\Lambda+N$.
Then, $s(M)<0$ if and only if $\rho(-\Lambda^{-1}N)<1$,
$s(M)=0$ if and only if $\rho(-\Lambda^{-1}N)=1$,
and $s(M)>0$ if and only if $\rho(-\Lambda^{-1}N)>1$.
\label{ff}
\end{proposition}



\begin{lemma}
{\rm (Proposition 2 in \cite{rantzer})}
Suppose that $M$ is a Metzler matrix such that $s(M)<0$.
Then, there exists a positive diagonal matrix $P$ such that
$M'P+PM$ is negative definite.
\label{less}\end{lemma}

\begin{lemma}
{\rm (Lemma A.1 in \cite{KhanaferAutomatica14})}
Suppose that $M$ is an irreducible Metzler matrix such that $s(M)=0$.
Then, there exists a positive diagonal matrix $P$ such that
$M'P+PM$ is negative semi-definite.
\label{equal}\end{lemma}

\section{The Bi-Virus Model}\label{system}


\begin{figure}
    \centering
    \begin{overpic}[width=\columnwidth]{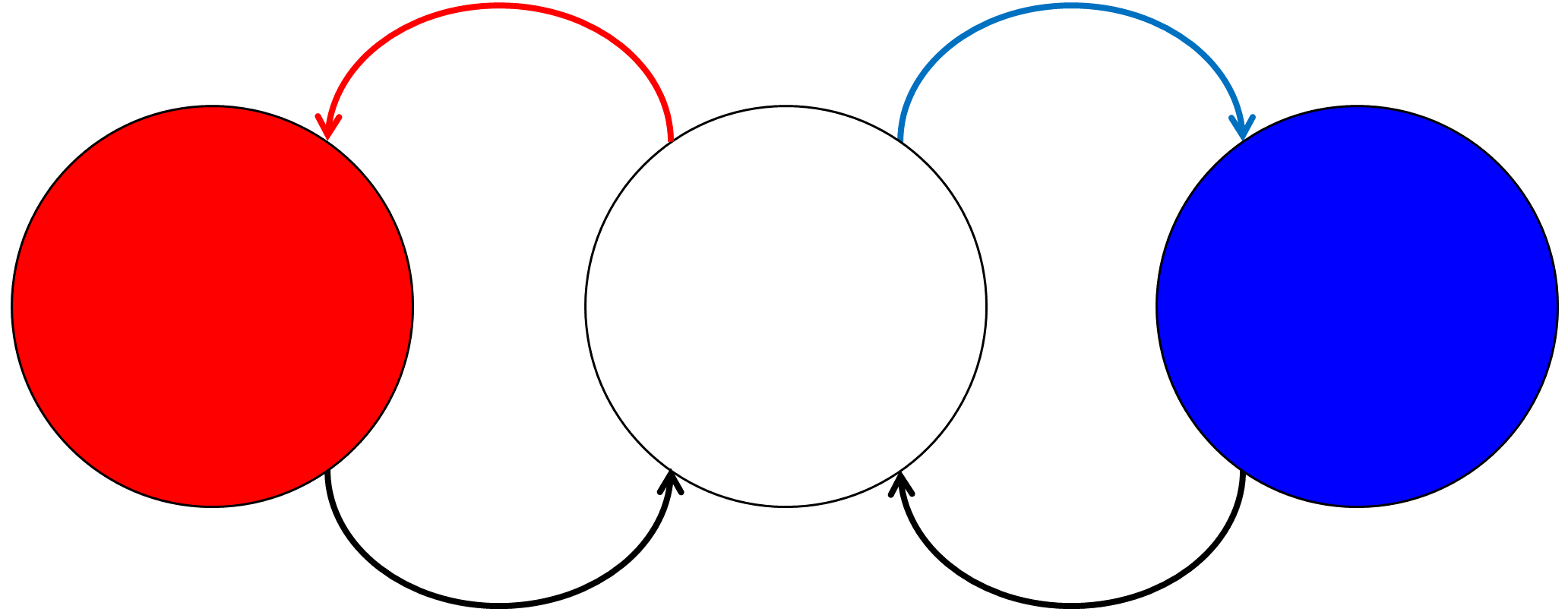}
        \put(30,5){{\parbox{0.05\linewidth}{%
     $$\delta_i^1$$}}}
     \put(-18,30){{\parbox{\linewidth}{%
     $$\sum \beta_{ij}^1  x^1_j$$}}}
     \put(18,30){{\parbox{\linewidth}{%
     $$\sum \beta_{ij}^2  x^2_j$$}}}
     \put(12.5,12){{\parbox{\linewidth}{%
     $x_i^1$}}}
     \put(12,20){{\parbox{\linewidth}{%
     {\color{white}\Huge I$^1$}}}}
     \put(47.75,19.2){{\parbox{\linewidth}{%
     \Huge S}}}
     \put(40.2,12){{\parbox{\linewidth}{%
     $1-x_i^1-x_i^2$}}}
     \put(85,20){{\parbox{\linewidth}{%
     {\color{white}\Huge I$^2$}}}}
     \put(85,12){{\parbox{\linewidth}{%
     $x_i^2$}}}
     \put(67,5){{\parbox{0.05\linewidth}{%
     $$\delta_i^2$$}}}
      \end{overpic}
    \caption{There are three states for each group $i$, with a portion of the group in each state:  a percentage of the group, $1 - x^{1}_i(t) - x^{2}_i(t)$, is susceptible (S); a percentage of the group,  $x^{1}_i(t)$, is infected with virus 1 (I$^{1}$); and a percentage of the group, $x^{2}_i(t)$, is infected with virus 2 (I$^{2}$). The healing and infection rates are indicated by $\delta_i^k$ and $\sum\beta_{ij}^kx_j^k$, respectively, for $k\in[2]$.}
    \label{fig:bivirus}
\end{figure}
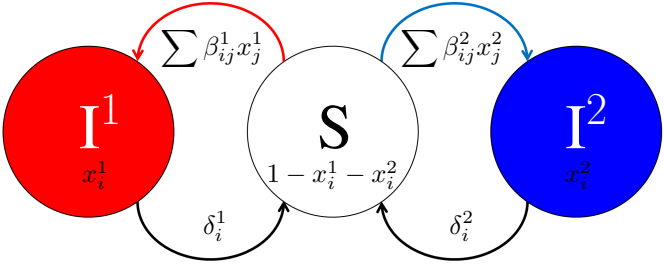

We are interested in the following continuous-time distributed model for two competing viruses.
Consider a network consisting of $n>1$ groups of individuals, labeled $1$ to $n$.
There are two competing viruses spreading over the network.
An individual cannot be infected with both viruses simultaneously.
An individual may be infected with one of the viruses, only by those in its own and neighboring groups.
Neighbor relationships among the $n$ groups are described
by a directed graph $\bbb{G}$ on $n$ vertices with an arc from vertex $j$
to vertex $i$ whenever the individuals in group $i$ can be infected by those in group $j$.
Thus, the neighbor graph $\bbb{G}$ has self-arcs at all $n$ vertices and the directions of arcs
in $\bbb{G}$ represent the directions of 
 contagion.
Each virus spreads over a spanning subgraph of $\bbb{G}$.
The two subgraphs can be different. Their union is the neighbor graph $\bbb G$.
It will be assumed that the two subgraphs are strongly connected
and, thus, so is $\bbb{G}$.\footnote{
A directed graph is {\em strongly connected} if
for every pair of distinct vertices $i$ and $j$, there is a directed path from $i$ to $j$ in the graph.
}

Let $S_i(t)$ denote the number of susceptible individuals in group $i$
at time $t\ge 0$, and
let $I^1_i(t)$ and $I^2_i(t)$ respectively denote the number of individuals infected
with virus 1 and virus 2 in group $i$ at time $t\ge 0$.
Assume that the total number of individuals in each group $i$, denoted by $N_i$, does not change over time.
In other words, $S_i(t)+I^1_i(t)+I^2_i(t)=N_i$, for all $i\in[n]$ and $t\ge 0$.
Several parameters are associated with each group $i$: healing rates $\gamma^1_i$ and $\gamma^2_i$
for virus 1 and virus 2 respectively,
birth rate $\mu_i$, death rate $\bar \mu_i$, and infection rates $\alpha^1_{ij}$ and $\alpha^2_{ij}$
for virus 1 and virus 2 respectively,
$i,j\in[n]$. 
Since $N_i$ is constant, $\bar\mu_i=\mu_i$.
We assume that individuals are susceptible at birth even if their parents are infected.
The evolution of the number of infected and susceptible individuals in each group $i$ is as follows:
\begin{align}\label{xxx1}
\dot S_i(t) &=  \mu_i N_i - \bar\mu_i S_i(t) + \gamma^1_iI^1_i(t) + \gamma^2_iI^2_i(t) \cr
         & \;\;\; - \sum_{j=1}^n \alpha^1_{ij} \frac{S_i(t)}{N_i} I^1_j(t) - \sum_{j=1}^n \alpha^2_{ij} \frac{S_i(t)}{N_i} I^2_j(t) \cr
         &=  (\mu_i + \gamma^1_i )I^1_i(t) +(\mu_i + \gamma^2_i )I^2_i(t) \cr
         & \;\;\; - \sum_{j=1}^n \alpha^1_{ij} \frac{S_i(t)}{N_i} I^1_j(t)  - \sum_{j=1}^n \alpha^2_{ij} \frac{S_i(t)}{N_i} I^2_j(t),
 \end{align}
\begin{align} \label{xxx2}
\dot I^1_i(t) &= - \gamma^1_iI^1_i(t) - \bar\mu_i I^1_i(t) + \sum_{j=1}^n \alpha^1_{ij} \frac{S_i(t)}{N_i} I^1_j(t) \cr
         &= (-\gamma^1_i-\mu_i) I^1_i(t) + \sum_{j=1}^n \alpha^1_{ij} \frac{S_i(t)}{N_i} I^1_j(t),
\end{align}
\begin{align} \label{xxx3}
\dot I^2_i(t) &= - \gamma^2_iI^2_i(t) - \bar\mu_i I^2_i(t) + \sum_{j=1}^n \alpha^2_{ij} \frac{S_i(t)}{N_i} I^2_j(t) \cr
         &= (-\gamma^2_i-\mu_i) I^2_i(t) + \sum_{j=1}^n \alpha^2_{ij} \frac{S_i(t)}{N_i} I^2_j(t),
\end{align}
where the infection of group $i$ is caused by one of its neighboring groups $j$, proportional to the total number of infected individuals in group $j$ and the proportion of susceptible individuals in group $i$, which can be regarded as the probability of contact. It is worth noting that, since group $i$ is defined as a neighboring group of itself, its infected individuals can infect its own susceptible members, which reflects realistic scenarios. It is easy to see that since $S_i(t)+I^1_i(t)+I^2_i(t)=N_i$ is a constant, \rep{xxx1} can be implied by \rep{xxx2} and \rep{xxx3}.
To simplify the model, define the proportion of infected individuals in group $i$ by
$$x^1_i(t) = \frac{I^1_i(t)}{N_i},  \;\;\;\;\; x^2_i(t) = \frac{I^2_i(t)}{N_i},$$
and let
$$\beta^1_{ij} = \alpha^1_{ij} \frac{N_j}{N_i}, \;\;\; \beta^2_{ij} = \alpha^2_{ij} \frac{N_j}{N_i}, \;\;\; \delta^1_i=\gamma^1_i+\mu_i, \;\;\; \delta^2_i=\gamma^2_i+\mu_i.$$
From \rep{xxx2} and \rep{xxx3}, it follows that
\begin{equation}\label{updates}
\begin{split}
\dot{x}^{1}_i(t) &= - \delta^{1}_i x^{1}_i(t) + (1 - x^{1}_i(t) - x^{2}_i(t))\sum_{j=1}^n \beta^{1}_{ij}  x^{1}_j(t) , \\
\dot{x}^{2}_i(t) &= - \delta^{2}_i x^{2}_i(t) +(1 - x^{2}_i(t) - x^{1}_i(t))\sum_{j=1}^n \beta^{2}_{ij}  x^{2}_j(t) .
\end{split}
\end{equation}
The progression from \rep{xxx3} to the second equation in \rep{updates} is given as follows. 
Dividing both sides of \rep{xxx3} by $N_i$, we have
$$\frac{\dot I^2_i(t)}{N_i} = (-\gamma^2_i-\mu_i) \frac{I^2_i(t)}{N_i} + \sum_{j=1}^n \alpha^2_{ij} \frac{S_i(t)}{N_i} \frac{N_j}{N_i}\frac{I^2_j(t)}{N_j}.$$
Since $S_i(t)+I^1_i(t)+I^2_i(t)=N_i$, we have 
$$\frac{\dot I^2_i(t)}{N_i} = (-\gamma^2_i-\mu_i) \frac{I^2_i(t)}{N_i} + \frac{N_i-I^1_i(t)-I^2_i(t)}{N_i}\sum_{j=1}^n \alpha^2_{ij}  \frac{N_j}{N_i}\frac{I^2_j(t)}{N_j}.$$
Since 
$$x^1_i(t) = \frac{I^1_i(t)}{N_i},  \;\;\; x^2_i(t) = \frac{I^2_i(t)}{N_i}, \;\;\; \beta^2_{ij} = \alpha^2_{ij} \frac{N_j}{N_i}, \;\;\; \delta^2_i=\gamma^2_i+\mu_i,$$
it follows that 
$$\dot{x}^{2}_i(t) = - \delta^{2}_i x^{2}_i(t) +(1 - x^{2}_i(t) - x^{1}_i(t))\sum_{j=1}^n \beta^{2}_{ij}  x^{2}_j(t),$$
which is the same as the second equation in \rep{updates}. The first equation in \rep{updates} can be derived from \rep{xxx2} in the same way.

The above derivation generalizes the one in \cite{FallMMNP07} for a single virus to two competing viruses. See \cite{FallMMNP07} for a more detailed explanation of the derivation.

Note that each virus has its own non-symmetric infection rates $\beta^{1}_{ij},\beta^{2}_{ij}$ incorporating the nearest-neighbor graph structures and healing rates $\delta^{1}_i,\delta^{2}_i$. 
A graphical depiction of this model is given in Figure \ref{fig:bivirus}. The model can be written in a matrix form as
\begin{equation}\label{sys}
\begin{split}
\dot{x}^{1}(t) &= (- D^{1}+B^{1} - X^{1}(t)B^{1} - X^{2}(t)B^{1} )x^{1}(t), \\
\dot{x}^{2}(t) &= (- D^{2}+B^{2} - X^{2}(t)B^{2} - X^{1}(t)B^{2} )x^{2}(t),
\end{split}
\end{equation}
where ${x}^{k}(t)\in [0,1]^n$, $B^{k}$ is the matrix of $\beta^{k}_{ij}$'s, $X^{k}(t)={\rm diag}(x^{k}(t))$, and $D^{k}={\rm diag}(\delta^{k})$, with $k=1,2$ indicating virus $1$ or $2$.

The same mathematical model was first proposed in \cite{prakash2012winner} with an alternative  interpretation
in which the system consists of $n$ agents, and
${x}^{1}_i(t)$ and ${x}^{2}_i(t)$ are the probabilities that agent $i$ has viruses $1$ and $2$, respectively.
The model can be viewed as a simplified model resulting from a mean field
approximation on a $3^n$ state Markov chain model, similar to what has been done for the single-virus SIS model in \cite{OmicTN09}.

\begin{figure}
    \centering
    \includegraphics[width=\columnwidth]{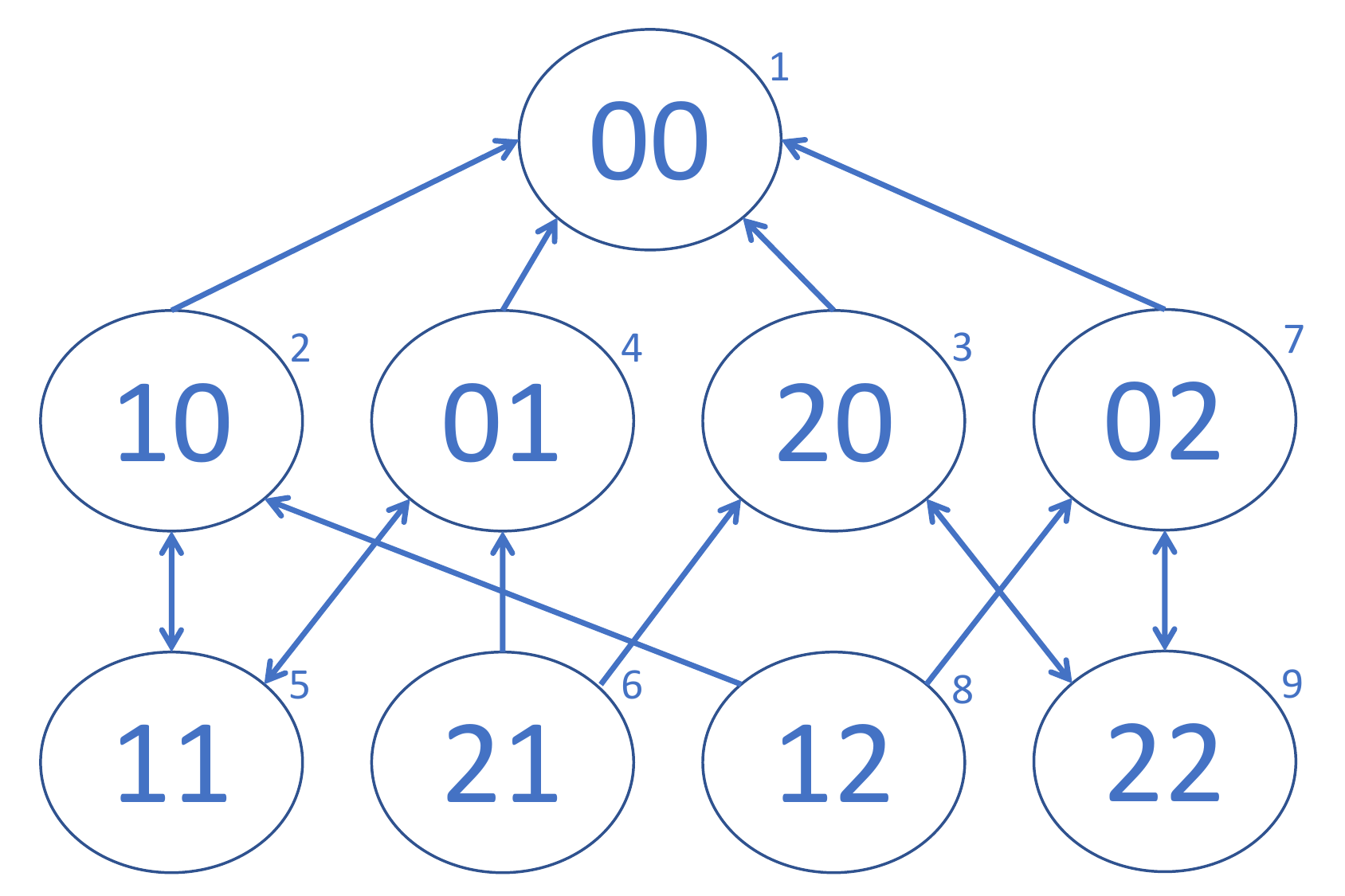}
    \caption{
    Example of $3^n$-state model with $n=2$: the superscripts indicate the ordering of the states, which correspond to the subscript of $y_k(t)$ in \eqref{eq:y}, and the internal strings indicate which agents are healthy (state 0) and which are infected with virus~1 (state 1) or virus 2 (state 2), corresponding to $s_i$ in \eqref{eq:qij} (the $i$th entry from left to right corresponds to the $i$th agent).}
    \label{fig:3n}
\end{figure}
For completeness, we provide here a full description of the $3^n$ state Markov model. 
Each state of the chain, $Y_k(t)$, corresponds to a ternary-valued string $s$ of length $n$, where  $s_i=0$, or $s_i=1$, or $s_i=2$ indicate that the $i$th agent is either susceptible, or infected with virus 1, or infected with virus 2, respectively. 
The state transition matrix, $Q$, is defined by
\begin{equation}\label{eq:qij}
{q}_{kl}=
\begin{cases}
\delta_i^1, & \text{ if } s_i = 1, k = l + 3^{i-1}\\
\delta_i^2, & \text{ if } s_i = 2, k = l + 2(3^{i-1})\\
\displaystyle\sum_{j=1}^n \beta^1_{ij}  1_{s_j=1}, &\text{ if } s_i = 0, k = l - 3^{i-1} \\[2.5ex]
\displaystyle\sum_{j=1}^n \beta^2_{ij}  1_{s_j=2}, &\text{ if } s_i = 0, k = l - 2(3^{i-1}) \\[2.5ex]
-\displaystyle\sum_{j\neq l} {q}_{jl}, & \text{ if } k = l\\
0, & \text{ otherwise,}
\end{cases}
\end{equation}
for $i\in [n]$. Here virus 1 and virus 2 are propagating over a network whose infection rates are given by $\beta^1_{ij}$ and $\beta^2_{ij}$, respectively (nonnegative with $\beta^1_{ii} = \beta^2_{ii}=0 ,\ \forall j$), 
 $\delta_i^1$ and $\delta_i^1$ are the respective healing rates of the $i$th agent, and, again, $s_i=0$, or $s_i=1$, or $s_i=2$ indicate that the $i$th agent is either susceptible, or infected with virus 1, or infected with virus 2, respectively. 
The state vector $y(t)$ is defined as 
\begin{equation}\label{eq:y}
         y_k(t) = Pr[Y_k(t) = k], 
\end{equation}
with $\sum_{k=1}^{3^n} y_k(t) =1$. The Markov chain evolves as
\begin{equation}\label{eq:3n}
    \frac{dy'(t)}{dt} = y'(t){Q}.
\end{equation}
See Figure \ref{fig:3n} for an illustration of this chain with $n=2$.


Let  $v^1_i(t) = Pr[X_i(t) = 1]$ and $v^2_i(t) = Pr[X_i(t) = 2]$, where $X_i(t)$ is the random variables representing whether the $i$th agent is susceptible or infected with virus 1 or 2. Then
\begin{equation}\label{eq:v}
\begin{split}
    (v^1)'(t) = y'(t)M^1,\\
	(v^2)'(t) = y'(t)M^2,
\end{split}
\end{equation}
where the $i$th columns of $M^1$ and $M^2$ indicate the states in the Markov chain where agent $i$ is infected with virus 1 and 2 (all the ternary strings where $s_i = 1$ and 
$s_i = 2$), respectively, that is,
\begin{align*}
M^1 &= 1_{M=1},\\
M^2 &= 1_{M=2},
\end{align*}
where $M\in \mathbb{R}^{3^n\times n}$ has rows of lexicographically-ordered ternary numbers, bit reversed.\footnote{Matlab code: $M = fliplr(dec2base(0:(3^n)-1,3)-'0')$} 
Therefore, $v^1_i(t)$ and $v^2_i(t)$ reflect the summation of all probabilities where $s_i = 1$ and $s_i = 2$. 
Note that the first state of the chain, which corresponds to $s = \0$, the healthy state, for $\delta^1_i,\delta^2_i>0 \ \forall i$, is the absorbing, or sink, state of the chain. This means that the Markov chain will never escape the state once in it, and further, since it is the only absorbing state the system will converge to the healthy state with probability one \cite{norris1998markov}. 

In \cite{FallMMNP07}, a traditional single-group deterministic SIS model was generalized to a networked multi-group setting, similar to the process in \eqref{xxx1}-\eqref{updates}, which results in the same model as the mean-field approximation of the networked Markov chain model proposed in \cite{OmicTN09}, except that $\beta_{ii}$ can be nonzero in the deterministic case, which is not possible in the probabilistic derivation, and therefore is more general. Similar to \cite{OmicTN09}, the model in \eqref{sys} can derived as a first order approximation of the $3^n$ state Markov chain model in \eqref{eq:qij}-\eqref{eq:v}, which is the model studied in \cite{watkins2016optimal}. Consequently, the states of \eqref{sys} can be interpreted as the probability of the agents being infected or the proportion of subpopulations that are infected. Since using the deterministic bi-virus model for a networked multi-group setting allows $\beta_{ii}$ to be nonzero and is more general, we focus here on this model. 
For completeness, to illustrate the effectiveness of the first order approximation, we compare \eqref{sys} and \eqref{eq:qij}-\eqref{eq:v} via simulations in Section \ref{simulations}. 


Note that if ${x}^{2}(t) = \0$ for all $t$, \eqref{updates} reduces to the single-virus model (we drop the superscript since there is only one virus),
\eq{\dot x_i(t)=-\delta_i x_i(t)+(1-x_i(t))\sum_{j=1}^n \beta_{ij}x_j(t),
\label{update}}
where $x_i(t)$ is the proportion of infected individuals
in group $i$ (or the probability that agent $i$ has the virus), 
$\beta_{ij}$'s are the infection rates, $\delta_i$'s are the healing rates, and $x_i(0)\in[0,1],i\in[n]$,
or in matrix form
\eq{\dot x(t) = \left(-D+B-X(t)B \right)x(t). \label{single}}
Consider a further special case in which all $n$ groups are isolated, i.e., $\beta_{ij}=0$ for all $i\neq j$. For each group $i$, if $\beta_{ii}>0$, from \rep{update}, 
$$\dot x_i(t)=-\delta_i x_i(t)+\beta_{ii}(1-x_i(t))x_i(t),$$
which is an SIS model for a single group of individuals. It is not hard to verify that for this single-group SIS model, if $\delta_i\geq \beta_{ii}$, $x_i(t)$ will always converge to zero, and if $\delta_i < \beta_{ii}$, $x_i(t)$ will always converge to a positive value unless $x_i(0)=0$. If $\beta_{ii}=0$, then the dynamics reduce to $\dot x_i(t)=-\delta_i x_i(t)$ which is not an epidemic model. Therefore, $\beta_{ii}>0$ has an important physical meaning in networked models of groups of individuals (sometimes called metapopulation models in the literature).
If we interpret each state of the model to be the infection probability of a single agent, and we further assume that $\beta_{ij}$ can be factored into $\beta_i a_{ij}$, where $\beta_i$ is the infection rate of agent $i$, and $a_{ij}$ is the connection structure between agents, then \rep{update} becomes the single SIS model proposed in~\cite{OmicTN09}.

We impose the following assumptions on the model
throughout Sections~\ref{equilibria} and \ref{sense}.

\begin{assumption}
For all $i\in[n]$, we have $x^{1}_i(0),x^{2}_i(0),(1-x^{1}_i(0)-x^{2}_i(0))\in[0,1]$.
\label{x0}
\end{assumption}

\begin{assumption}
For all $i\in[n]$, we have $\delta^{1}_i,\delta^{2}_i\ge0$. The matrices $B^{1}$ and $B^{2}$ are nonnegative and irreducible.
\label{para}
\end{assumption}

Assumption~\ref{x0} says that the initial proportions of infected and healthy individuals are in the interval $[0,1]$.
The nonnegativity assumption on the  matrix $B^{k}$ is equivalent to
$\beta^{k}_{ij}\ge 0$ for all $k\in[2]$ and $i,j\in[n]$.
Assumption~\ref{para} says that all healing and infection rates are nonnegative.
The assumption of an irreducible matrix $B^{k}$ is equivalent to a strongly connected
spreading graph for virus $k$, $k\in[2]$.


\begin{lemma}
Under the conditions of Assumptions  \ref{x0} and \ref{para}, $x^{1}_i(t),x^{2}_i(t),x^{1}_i(t)+x^{2}_i(t)\in[0,1]$ for all $i\in[n]$ and $t\ge 0$.
\label{box}
\end{lemma}

Lemma \ref{box} implies that the set
\begin{eqnarray}\label{D}
    \D 
    =\{(x^1, x^2) \; | \; x^1\ge \0, \; x^2\ge \0, \; x^1+x^2\le \1\}
\end{eqnarray}
is positively invariant with respect to the system defined by \rep{sys}. Since $x^{1}_i$ and $x^{2}_i$
denote the fractions  of group $i$ infected by viruses 1 and 2, respectively,
and $1-x^{1}_i-x^{2}_i$ denotes the fraction of group $i$ that is healthy, it is natural to assume that their initial values are in the interval $[0,1]$,
since otherwise the values will lack any physical meaning for the epidemic model considered here.
Therefore, in this paper,
we focus on the analysis of \rep{sys} only on the domain $\D$, as defined in \rep{D}.

We are interested in the problem of characterizing limiting behavior of the bi-virus model \rep{sys} and its dependence on the network structure (two spreading graphs) and parameters (healing rates $\delta_i^1$, $\delta_i^2$ and infection rates $\beta_{ij}^1$, $\beta_{ij}^2$). 
The limiting behaviors will be characterized by the equilibria of the system and their stability. The effects of the network structure and parameters on the limiting behavior are important and useful for controlling the epidemic spreading process. Specifically, it will be shown in the following sections that, under certain assumptions, whether the two viruses 
eventually disappear or not can be determined by checking the network structure and parameters. In the case when at least one virus 
ultimately spreads over the network, the spreading process can be attenuated or eliminated by modifying the network parameters.

\section{Equilibria and Their Stability} \label{equilibria}

In this section, we analyze the equilibria of the system \rep{sys} and their stability, which characterize limiting behaviors of the bi-virus model.

First, it can be seen that $x^{1}=x^{2}=\0$ is an equilibrium of the system \rep{sys},
which corresponds to the case when  no individual is infected.
We call this trivial equilibrium the {\em healthy state}.
We will show that \rep{sys} also admits nonzero equilibria under appropriate assumptions.
We call those nonzero equilibria {\em epidemic states}.
In this section, we study the stability of the healthy state as well as the epidemic states of \rep{sys}.
To state our results, we need the following definition.

\begin{definition}
Consider an autonomous system
\eq{\dot x(t)  = f(x(t)), \label{def}}
where $f: \scr{X}\rightarrow\R^n$ is a locally Lipschitz map from a domain $\scr{X}\subset\R^n$ into $\R^n$.
Let $z$ be an equilibrium of \rep{def} and $\scr{E}\subset\scr{X}$ be a domain containing $z$.
If the equilibrium $z$ is asymptotically stable such that for any $x(0)\in\scr{E}$ we have
$\lim_{t\rightarrow\infty}x(t) = z$, then $\scr{E}$ is said to be a domain of attraction for $z$.
\end{definition}

The following corollary is a direct consequence of Lyapunov's stability theorem
(see Theorem 4.1 in \cite{khalil}) and the definition of domain of attraction.

\begin{corollary}
Let $z$ be an equilibrium of \rep{def} and $\scr{E}\subset \scr{X}$ be a domain
containing $z$. Let $V:\scr{E}\rightarrow\R$ be a continuously differentiable function
such that $V(z)=0$, $V(x)>0$ in $\scr{E}\setminus \{z\}$, $\dot V(z)=0$,
and $\dot V(x)<0$ in $\scr{E}\setminus \{z\}$. If $\scr{E}$ is a positively invariant set,
then the equilibrium $z$ is asymptotically stable with a domain of attraction $\scr{E}$.
\label{lya}
\end{corollary}

The following theorem establishes a sufficient condition for global stability of the healthy state, whereas the work of \cite{watkins2015optimal,watkins2016optimal} provides 
necessary and sufficient conditions for local stability of the healthy state.

\begin{theorem}
Under Assumption \ref{para}, 
the healthy state is the unique equilibrium of system~\rep{sys} if, and only if,
$s(-D^{1}+B^{1})\leq 0$ and $s(-D^{2}+B^{2})\leq 0$. Furthermore, in this case,
the healthy state is
asymptotically stable with  domain of attraction $\D$, as defined in \rep{D}.
\label{0global}\end{theorem}

To prove the theorem, we need the following 
result
for the single-virus model \rep{single}.

\begin{proposition}
Suppose that $\delta_i\ge0$ for all $i\in[n]$ and that matrix $B$ is nonnegative and irreducible.
If $s(-D+B)\leq 0$, then
$\0$ is the unique equilibrium of system~\rep{single}, which is asymptotically stable with domain of attraction $[0,1]^n$.
If $s(-D+B)> 0$, then system~\rep{single} has two equilibria, $\0$ and $x^*$ which satisfies $x^*\gg \0$.
\label{single0}\end{proposition}

This result has been proved in \cite{Gonorrhea,FallMMNP07,KhanaferAutomatica14} for the case
when $\delta_i > 0$ for all $i\in[n]$. Additional assumptions on $\beta_{ij}$ are also imposed in \cite{Gonorrhea,KhanaferAutomatica14}. Using a different Lyapunov function from the one used in the aforementioned papers, we extend the result by allowing $\delta_i=0$ (see the Appendix), which reflects the situations where certain groups of individuals are unable to heal themselves (due, for example, to the lack of vaccines). 
Note that this allows for the SI model in which $\delta_i = 0$ for all $i\in[n]$.
In the case when $\delta_i > 0$ for all $i\in[n]$, the conditions $s(-D + B) \le 0$ and $\rho(D^{-1}B) \le 1$ are equivalent, which is a direct consequence of Proposition~\ref{ff}. 
In the case when $\delta_i = 0$ for some but not all $i\in[n]$, the two conditions are not equivalent since $\rho(D^{-1}B)$ is not well-posed in this case.






{\em Proof of Theorem \ref{0global}:}
We first show that
both $x^{1}(t)$ and $x^{2}(t)$ will asymptotically converge to $\0$ as $t\rightarrow\infty$
for any initial condition in $\D$.
Since $x^{1}_i(t)$ and $x^{2}_i(t)$ are always nonnegative by Lemma \ref{box},
from \rep{updates}, we obtain
\begin{eqnarray*}
\dot{x}^{1}_i(t) &\le& - \delta^{1}_i x^{1}_i(t) + (1 - x^{1}_i(t))\sum_{j=1}^n \beta^{1}_{ij}  x^{1}_j(t) , \\
\dot{x}^{2}_i(t) &\le& - \delta^{2}_i x^{2}_i(t) + (1 - x^{2}_i(t))\sum_{j=1}^n \beta^{2}_{ij}  x^{2}_j(t) ,
\end{eqnarray*}
which imply that each of the trajectories of $x^{1}_i(t)$ and $x^{2}_i(t)$ is bounded above by
a single-virus model. From Assumption \ref{para} and Proposition \ref{single0},
both $x^{1}_i(t)$ and $x^{2}_i(t)$ will asymptotically converge to $\0$ as $t\rightarrow\infty$,
and thus the healthy state is the unique equilibrium of \rep{sys}.

We next show the asymptotic stability of the healthy state. 
Consider the Lyapunov function candidate
$V(x^1(t),x^2(t)) = x^1(t)'P^1x^1(t) + x^2(t)'P^2x^2(t)$, where $P^1$ and $P^2$ are positive diagonal matrices chosen in the same way as $P$ is chosen in Proposition \ref{single0}. 
Note from \eqref{sys} that
\begin{eqnarray*}
&&\dot V(x^1(t),x^2(t)) \\
&=& 2x^1(t)'P^1(-D^1+B^1-X^1(t)B^1-X^2(t)B^1)x^1(t)\\
&&+2x^2(t)'P^2(-D^2+B^2-X^1(t)B^2-X^2(t)B^2)x^2(t)\\
&=& x^1(t)'((-D^1+B^1)'P^1+P^1(-D^1+B^1))x^1(t) \\
&& +x^2(t)'((-D^2+B^2)'P^2+P^2(-D^2+B^2))x^2(t) \\
&& -2x^1(t)'P^1(X^1(t)B^1+X^2(t)B^1)x^1(t) \\
&& -2x^2(t)'P^2(X^1(t)B^2+X^2(t)B^2)x^2(t) \\
&\le& x^1(t)'((-D^1+B^1)'P^1+P^1(-D^1+B^1))x^1(t) \\
&& +x^2(t)'((-D^2+B^2)'P^2+P^2(-D^2+B^2))x^2(t).
\end{eqnarray*}
Using similar arguments to those in the proof of Proposition \ref{single0}, it can be shown that $\dot V(x^1(t),x^2(t))$ is a negative definite function over $\scr D$ except for the healthy state. 
By Lemma \ref{box} and Corollary \ref{lya},
the healthy state is asymptotically stable with  domain of attraction $\scr D$.
Now we show that if either $s(-D^{1}+B^{1})> 0$ or $s(-D^{2}+B^{2})> 0$,
then system \rep{sys} has an epidemic state.

Without loss of generality, suppose that $s(-D^{1}+B^{1})> 0$. Set $x^{2}=\0$. Then,
the dynamics of $x^{1}$ simplifies to that of the single-virus system, which admits an epidemic state
by Proposition~\ref{single0}. Therefore, in the case when  $s(-D^{1}+B^{1})> 0$,
the system \rep{sys} always admits an equilibrium of the form $(\tilde x^{1}, \0)$ with $\tilde x^{1}\gg \0$.
\hfill
$\qed$


We have provided a necessary and sufficient condition for the eradication of both of the viruses,
i.e., $s(-D^{1}+B^{1})\le 0$ and $s(-D^{2}+B^{2})\leq 0$.
Since larger nonzero entries of $D^1$ and $D^2$ (i.e., healing rates $\delta_i^1$, $\delta_i^2$) will decrease the two quantities, and larger nonzero entries of $B^1$ and $B^2$ (i.e., infection rates $\beta_{ij}^1$, $\beta_{ij}^2$) will increase the two quantities, the necessary and sufficient condition can be interpreted as 
the overall healing capabilities of all the individuals overcoming or balancing out competely the effects of the network infection.
This characterization is important for understanding when a system will become completely healthy as illustrated via simulation in Figure \ref{fig:dfe}.

Now we turn to the analysis of epidemic states. We begin with dominant virus states at which one
virus is eradicated and the other one pervades in the network.

\begin{theorem}
Under Assumption \ref{para}, if $s(-D^{1}+B^{1})> 0$ and $s(-D^{2}+B^{2})\leq 0$, then
\rep{sys} has two equilibria,
the healthy state $(\0,\0)$, 
where the system converges to this equilibrium for all initial conditions $(x^{1}(0),x^{2}(0)) \in\{(\0,x^2) | x^2 \in [0,1]^n\}$, 
and a unique epidemic state of the form $(\tilde x^{1}, \0)$ with $\tilde x^{1}\gg \0$, which is
asymptotically stable with  domain of attraction $\D\setminus\{(\0,x^2) | x^2 \in [0,1]^n\}$, with $\D$ defined in \rep{D}.
\label{eglobal}\end{theorem}

\begin{remark}
Note that healthy state $(\0,\0)$ is an unstable equilibrium. A small perturbation of the first virus from the origin will drive the system to the unique epidemic state $(\tilde x^{1}, \0)$.
\hfill $\Box$
\end{remark}

To prove Theorem \ref{eglobal}, we need the following result for the single-virus model \rep{single}, which builds on Proposition \ref{single0}.

\begin{proposition}
Consider the single-virus model \rep{single}.
Suppose that $\delta_i\ge0$ for all $i\in[n]$, and that the matrix $B$ is nonnegative and irreducible.
If $s(-D+B)> 0$, then
the epidemic state $x^*\gg \0$ is asymptotically stable with  domain of attraction $[0,1]^n\setminus\{\0\}$.
\label{2local}\end{proposition}

This result has been proved in \cite{FallMMNP07,KhanaferAutomatica14} for the case
when $\delta_i > 0$ for all $i\in[n]$. We extend the result by allowing $\delta_i=0$.
The analyses in \cite{FallMMNP07,KhanaferAutomatica14} cannot be applied here.

To prove Proposition \ref{2local}, we need the following lemma.

\begin{lemma}
Consider the single-virus model \rep{single}.
Suppose that $\delta_i\ge0$ for all $i\in[n]$, and that the matrix $B$ is nonnegative and irreducible.
If $x(0)> \0$, then there exists a $\tau \ge 0$ such that $x(\tau)\gg \0$.
\label{attract}
\end{lemma}

Proofs of these results are provided in the appendix.

We are now in a position to provide the proof for Theorem~\ref{eglobal}.

{\em Proof of Theorem \ref{eglobal}:}
 By Proposition~\ref{single0}, $x^{2}(t)$ will asymptotically converge to $\0$
as $t\rightarrow\infty$ for all initial values $(x^{1}(0),x^{2}(0)) \in \{(\0,x^2) | x^2 \in [0,1]^n\}$, since the system reduces to the single virus case for that set of initial conditions.

From \rep{sys}, we have
$$\dot{x}^{1}(t) = (- D^{1}+B^{1} - X^{1}(t)B^{1}) x^{1}(t) - X^{2}(t)B^{1}x^{1}(t).$$
Thus, we can regard the dynamics of $x^{1}(t)$ as an autonomous system
\eq{\dot{x}^{1}(t) = (- D^{1}+B^{1} - X^{1}(t)B^{1}) x^{1}(t),\label{temp}}
with a vanishing perturbation $- X^{2}(t)B^{1}x^{1}(t)$, which converges to $\0$  as $t\rightarrow\infty$.
From Proposition \ref{2local}, the autonomous system \rep{temp} will asymptotically converge to
a unique epidemic state $\tilde x^{1} \gg \0$ for any $x^{1}(0)\in[0,1]^n\setminus\{\0\}$.

Let $y^1_i(t)=x^1_i(t)-\tilde x^1_i$ for all $i\in[n]$, or equivalently,
$y^1(t)=x^1(t)-\tilde x^1$.
Then,
\begin{eqnarray*}
\dot y^1(t) &=& \left( -D^1 + (I-\tilde X^1)B^1 -{\rm diag}(B^1x^1(t)) \right)y^1(t) \\
&& -X^2(t)B^1(y^1(t)+\tilde x^1).
\end{eqnarray*}
Let $f(y^1(t))=( -D^1 + (I-\tilde X^1)B^1 -{\rm diag}(B^1x^1(t)) )y^1(t)$ and
$g(t,y^1(t)) = -X^2(t)B^1(y^1(t)+\tilde x^1)$.
Consider the Lyapunov function candidate
$$V(y^1(t)) = \max_{k\in[n]}\frac{|y^1_k(t)|}{\tilde x^1_k}.$$
Then,
$\dot V(y^1(t)) = \frac{\partial V}{\partial y^1}f(y^1(t))
+\frac{\partial V}{\partial y^1}g(t,y^1(t))$.
From the proof of Proposition \ref{2local},
$\frac{\partial V}{\partial t} + \frac{\partial V}{\partial y^1}f(t,y^1)<0$
unless $y^1(t)=\0$ (i.e., $x^1(t)=\tilde x^1$).
Since $x^2(t)$ asymptotically converges to $\0$,
so does $\frac{\partial V}{\partial y^1}g(t,y^1(t))$.
This implies that after a sufficiently long time, $\dot V(y^1(t))<0$ if
$x^1(t)$ does not equal $\tilde x^1$.
Using the same argument as in the proof of Proposition \ref{2local},
$(x^{1}(t),x^{2}(t))$ will asymptotically converge to
the unique epidemic state $(\tilde x^{1},\0)$ for any $(x^{1}(0),x^{2}(0))\in\D\setminus\{(\0,x^2) | x^2 \in [0,1]^n\}$, with $\D$ defined in \rep{D}.
\hfill
$\qed$

Theorem \ref{eglobal} provides conditions under which one virus will pervade the network, and the other one will be driven out. Understanding this condition is useful for characterizing when a designer (marketer, politician, etc.) will consistently dominate a competitor. This behavior is illustrated in Figure~\ref{fig:1dom}.

It is clear from the preceding results that as long as one of $s(-D^{k}+B^{k})$, $k\in\{1,2\}$,
is less than or equal to zero, at most one virus will ultimately spread over the network.
A natural question is thus whether the two viruses can coexist when  $s(-D^{k}+B^{k})$, $k\in\{1,2\}$,
are both larger than zero.
In the following, we will partially answer this question.
We begin with a result regarding non-coexisting equilibria.

Let $(\tilde x^{1},\tilde x^{2})$ be an equilibrium of \rep{sys}.
Here, both $\tilde x^{1}$ and $\tilde x^{2}$ can be $\0$.
Then,
the Jacobian matrix of the equilibrium, denoted by $J(\tilde x^{1},\tilde x^{2})$, with $\tilde B^{i} = {\rm diag}(B^{i}\tilde x^{i})$, $i\in[2]$, is
\begin{align}\label{jacob}
&J(\tilde x^{1},\tilde x^{2}) = \\
&\footnotesize\begin{bmatrix}
(I - \tilde X^{1} -\tilde X^{2})B^{1} -D^{1} - \tilde B^{1} & -\tilde B^{1} \cr
- \tilde B^{2} & (I - \tilde X^{1} -\tilde X^{2})B^{2} -D^{2} - \tilde B^{2}\end{bmatrix}\normalsize. \nonumber
\end{align}

\begin{theorem}
Under Assumption \ref{para}, if $s(-D^{1}+B^{1})> 0$ and $s(-D^{2}+B^{2})>0$, then
\rep{sys} has at least three equilibria,
the healthy state $(\0,\0)$, and two epidemic states of the form $(\tilde x^{1}, \0)$ with $\tilde x^{1}\gg \0$
and $(\0,\tilde x^{2})$ with $\tilde x^{2}\gg \0$.
The healthy state $(\0,\0)$ is  unstable.
\label{noncoexist}\end{theorem}

{\em Proof:}
The existence of the two epidemic states is an immediate consequence of Proposition \ref{single0}.
The healthy state $(\0,\0)$ is always an equilibrium of \eqref{sys}.
Since by \eqref{jacob}  
$$
J(\0,\0) = \begin{bmatrix}
-D^{1}+B^{1} & 0 \cr
0 & -D^{2}+B^{2}\end{bmatrix},
$$
which is unstable as $s(-D^{1}+B^{1})> 0$ and $s(-D^{2}+B^{2})>0$,
the healthy state $(\0,\0)$ is unstable.
\hfill
$\qed$

Analysis of any of the other possible equilibria is challenging for general cases. In the sequel, we will thus consider two special cases, characterized by 
some additional assumptions on the underlying graphs for the spread of viruses. 
It turns out that
if for both $k \in \{1, 2 \}$, $s(-D^k + B^k)$ are 
larger than zero, then the existence of coexisting equilibria is not guaranteed, as shown in the following special case, where two homogeneous viruses spread on the same graph.

\begin{assumption}
Viruses 1 and 2 spread over the same strongly connected directed graph $\bbb{G}=([n],\scr{E})$, with $\delta^{1}_i=\delta^{1}>0$ and $\delta^{2}_i=\delta^{2}>0$ for all $i\in[n]$,
and $\beta^{1}_{ij}=\beta^{1}>0$ and $\beta^{2}_{ij}=\beta^{2}>0$ for all $i=j\in[n]$ and $(i,j)\in\scr{E}$.
\label{special1}\end{assumption}

Under Assumption \ref{special1}, it should be clear that
$D^{1}=\delta^{1}I$, $D^{2}=\delta^{2}I$, $B^{1}=\beta^{1}A$,
and $B^{2}=\beta^{2}A$,
where $A$ is the adjacency matrix of $\bbb{G}$, which is an irreducible Metzler matrix.

\begin{theorem}
Suppose that Assumptions \ref{x0} and \ref{special1} hold.
Then, coexisting equilibria exist only if
$\frac{\delta^{1}}{\beta^{1}}= \frac{\delta^{2}}{\beta^{2}}$.
\label{homo1}\end{theorem}

This result has been proved in \cite{prakash2012winner} for the case
when $\bbb{G}$ is an undirected graph. We extend the result by allowing $\bbb{G}$ to be directed.
To prove the theorem, we need the following lemma.

\begin{lemma}
Suppose that Assumptions \ref{x0} and \ref{special1} hold.
If $(\tilde x^{1},\tilde x^{2})$ is an equilibrium of \rep{sys},
then $\tilde x^{1}+\tilde x^{2}\ll \1$.
\label{nonzero}\end{lemma}

{\em Proof of Theorem \ref{homo1}:}
To prove the theorem, suppose that, to the contrary, there exists an equilibrium
$(\tilde x^{1},\tilde x^{2})$ such that $\tilde x^{1},\tilde x^{2}>\0$
in the case when $\frac{\delta^{1}}{\beta^{1}} \ne \frac{\delta^{2}}{\beta^{2}}$.
From \rep{sys} and Assumption \ref{special1},
\begin{eqnarray*}
(I-\tilde X^{1}-\tilde X^{2})A\tilde x^{1} = \frac{\delta^{1}}{\beta^{1}}\tilde x^{1}, \\
(I-\tilde X^{1}-\tilde X^{2})A\tilde x^{2} = \frac{\delta^{2}}{\beta^{2}}\tilde x^{2}.
\end{eqnarray*}
From Lemma \ref{nonzero}, $(I-\tilde X^{1}-\tilde X^{2})$ is a positive diagonal matrix,
and thus $(I-\tilde X^{1}-\tilde X^{2})A$ is also an irreducible Metzler matrix.
Since $\tilde x^{1},\tilde x^{2}>\0$, from Lemma \ref{metzler},
$s((I-\tilde X^{1}-\tilde X^{2})A)=\frac{\delta^{1}}{\beta^{1}} = \frac{\delta^{2}}{\beta^{2}}$,
which contradicts the hypothesis that $\frac{\delta^{1}}{\beta^{1}} \ne \frac{\delta^{2}}{\beta^{2}}$.
Therefore, coexisting equilibria may exist only if
$\frac{\delta^{1}}{\beta^{1}}= \frac{\delta^{2}}{\beta^{2}}$.
\hfill
$\qed$

For the following, without loss of generality, we assume $\frac{\delta^{1}}{\beta^{1}} > \frac{\delta^{2}}{\beta^{2}}$.

\begin{theorem}
Suppose that Assumptions \ref{x0} and \ref{special1} hold and that
$s(A)>\frac{\delta^{1}}{\beta^{1}} > \frac{\delta^{2}}{\beta^{2}}$.
Then, system \rep{sys} has three equilibria, the healthy state
$(\0,\0)$ which is unstable, $(\tilde x^{1},\0)$ with $\tilde x^{1}\gg \0$ which is unstable,
and $(\0,\tilde x^{2})$ with $\tilde x^{2}\gg \0$ which is locally exponentially stable.
\label{homo2}\end{theorem}

This result has been proved in \cite{prakash2012winner} for the case
when $\bbb{G}$ is an undirected graph. We extend the result by allowing $\bbb{G}$ to be directed with a  proof technique similar to \cite{prakash2012winner}.
In \cite{santos15}, a sufficient condition is established for the case in which $\delta^1=\delta^2=1$ and $\beta_{ij}^k$, $k\in\{1,2\}$, are heterogeneous (see Corollary~4 in \cite{santos15}).

{\em Proof:}
From Theorem \ref{homo1}, the system \rep{sys} cannot have any equilibria of the form
$(\tilde x^{1},\tilde x^{2})$ with $\tilde x^{1},\tilde x^{2}> \0$.
Thus, if $(\tilde x^{1},\tilde x^{2})$ is an equilibrium of \rep{sys}, at least one of
$\tilde x^{1}$ and $\tilde x^{2}$ equals $\0$.
It is clear that $(\0,\0)$ is always an equilibrium.
Suppose that $\tilde x^{1}=\0$ and $\tilde x^{2}> \0$. Then, from Proposition \ref{single0},
$\tilde x^{2}\gg \0$ and is unique.
Similarly, when $\tilde x^{1}>\0$ and $\tilde x^{2}= \0$,
$\tilde x^{1}\gg \0$ and is unique.
Thus, the system \rep{sys} has exactly three equilibria.

Next we turn to the  stability of the three equilibria.
Note that from Assumption \ref{special1}, the hypothesis
$s(A)>\frac{\delta^{1}}{\beta^{1}} > \frac{\delta^{2}}{\beta^{2}}$ implies
 $s(-D^{1}+B^{1}), s(-D^{2}+B^{2})>0$.
Then, from Theorem \ref{noncoexist}, the healthy state $(\0,\0)$ is  unstable.

From \rep{jacob}, the Jacobian at $(\tilde x^{1},\0)$ equals
$$\begin{bmatrix}
\beta^{1}(I - \tilde X^{1})A -\delta^{1}I - \beta^{1}{\rm diag}(A\tilde x^{1}) & -\beta^{1}{\rm diag}(A\tilde x^{1}) \cr
0 & \beta^{2}(I - \tilde X^{1} )A -\delta^{2}I \end{bmatrix}.
$$
From \rep{sys} and Assumption \ref{special1}, 
$(I - \tilde X^{1} )A\tilde x^{1} = \frac{\delta^{1}}{\beta^{1}}\tilde x^{1}$.
It follows from Lemma \ref{nonzero} that $(I - \tilde X^{1} )A$ is an irreducible Metzler matrix.
Then, from Lemma \ref{metzler}, 
$s((I - \tilde X^{1} )A)= \frac{\delta^{1}}{\beta^{1}}$.
Since $\frac{\delta^{1}}{\beta^{1}} > \frac{\delta^{2}}{\beta^{2}}$,
it follows that
\begin{eqnarray*}
s(\beta^{2}(I - \tilde X^{1} )A -\delta^{2}I) &=& \beta^{2}s((I - \tilde X^{1} )A) -\delta^{2} \\
&= & \beta^{2}\left(\frac{\delta^{1}}{\beta^{1}} - \frac{\delta^{2}}{\beta^{2}}\right) > 0,
\end{eqnarray*}
which implies that the Jacobian matrix is unstable. Thus,
the equilibrium $(\tilde x^{1},\0)$ with $\tilde x^{1}\gg \0$ is  unstable.

From \rep{jacob}, the Jacobian at $(\0,\tilde x^{2})$ equals
$$\begin{bmatrix}
\beta^{1}(I - \tilde X^{2} )A -\delta^{1}I  &  0 \cr
-\beta^{2}{\rm diag}(A\tilde x^{2}) & \beta^{2}(I - \tilde X^{2})A -\delta^{2}I - \beta^{2}{\rm diag}(A\tilde x^{2}) \end{bmatrix}.
$$
Using the same arguments as in the previous paragraph,
 $s(\beta^{1}(I - \tilde X^{2} )A -\delta^{1}I)<0$.
From \rep{sys} and Assumption \ref{special1}, 
$$(I - \tilde X^{2} )A\tilde x^{2} = \frac{\delta^{2}}{\beta^{2}}\tilde x^{2}.$$
Since $\tilde x^{2}\gg \0$ and $A$ is irreducible,
it must be true that
$$\left(\beta^{2}(I - \tilde X^{2})A -\delta^{2}I - \beta^{2}{\rm diag}(A\tilde x^{2})\right)\tilde x^{2}<\0.$$
It follows from Assumptions \ref{x0} and \ref{special1} and Lemma \ref{box} 
that $\beta^{2}(I - \tilde X^{2})A -\delta^{2}I - \beta^{2}{\rm diag}(A\tilde x^{2})$ is an irreducible Metzler matrix.
Then, from Lemma \ref{metzler}, 
$s(\beta^{2}(I - \tilde X^{2})A -\delta^{2}I - \beta^{2}{\rm diag}(A\tilde x^{2}))<0$,
which implies that the Jacobian matrix is stable. Thus,
the equilibrium $(\0,\tilde x^{2})$ with $\tilde x^{2}\gg \0$ is locally exponentially stable.
\hfill
$\qed$

For the possibility of coexisting equilibria, we have the following interesting result.

\begin{theorem}
Suppose that Assumptions \ref{x0} and \ref{special1} hold and that $s(A)>\frac{\delta^{1}}{\beta^{1}} = \frac{\delta^{2}}{\beta^{2}}$.
If $(\tilde x^{1}, \tilde x^{2})$ with $\tilde x^{1},\tilde x^{2}>\0$ is an equilibrium
of \rep{sys}, then $\tilde x^{1},\tilde x^{2}\gg\0$ and
$\tilde x^{1} = \alpha\tilde x^{2}$ for some constant $\alpha>0$.
Furthermore, for each $\alpha>0$ there exists a unique pair $(\tilde x^{1}, \tilde x^{2})$ such that $\tilde x^{1} = \alpha\tilde x^{2}$.
\label{parallel1}\end{theorem}

{\em Proof:}
From the proof of Theorem \ref{homo1}, 
\begin{align}\label{equil}
\begin{split}
(I-\tilde X^{1}-\tilde X^{2})A\tilde x^{1} = \frac{\delta^{1}}{\beta^{1}}\tilde x^{1}, \\
(I-\tilde X^{1}-\tilde X^{2})A\tilde x^{2} = \frac{\delta^{2}}{\beta^{2}}\tilde x^{2},
\end{split}
\end{align}
in which $(I-\tilde X^{1}-\tilde X^{2})A$ is an irreducible Metzler matrix.
From Lemma \ref{metzler0}, it must be true that $\tilde x^{1},\tilde x^{2}\gg\0$ and
$\tilde x^{1} = \alpha\tilde x^{2}$ for some constant $\alpha>0$.

Given some $\alpha>0$, assume $\tilde x^{1} = \alpha\tilde x^{2}$ and $\hat x^{1} = \alpha\hat x^{2}$ both satisfy \eqref{equil}.
Therefore,
\begin{align*}
\begin{split}
\left( A -(1+\frac{1}{\alpha})\tilde X^{1}A\right)\tilde x^{1} = \frac{\delta^{1}}{\beta^{1}}\tilde x^{1}, \\
\left( A -(1+\frac{1}{\alpha})\hat   X^{1}A\right)\hat   x^{1} = \frac{\delta^{1}}{\beta^{1}}\hat   x^{1},
\end{split}
\end{align*}
so that, by Lemma \ref{metzler},
\begin{equation}\label{eq:s}
    s\left(A -(1+\frac{1}{\alpha})\tilde X^{1}A\right) = s\left(A -(1+\frac{1}{\alpha})\hat X^{1}A\right) = \frac{\delta^{1}}{\beta^{1}}.
\end{equation}
Also, without loss of generality, assume there exists $j \in [n]$ such that $\tilde x^{1}_j >  \hat x^{1}_j$ and $\tilde x^{1}_i =  \hat x^{1}_i$ for all $i\neq j$.
This implies
$$A -(1+\frac{1}{\alpha})\tilde X^{1}A < A -(1+\frac{1}{\alpha})\hat X^{1}A,$$
which by Lemma \ref{metzler} implies
$$s\left(A -(1+\frac{1}{\alpha})\tilde X^{1}A\right) < s\left(A -(1+\frac{1}{\alpha})\hat X^{1}A\right).$$
However, this contradicts \eqref{eq:s}.
So for each $\alpha>0$ there exists a unique pair $(\tilde x^{1}, \tilde x^{2})$ such that $\tilde x^{1} = \alpha\tilde x^{2}$.
\hfill
$\qed$

\begin{remark}\label{rem:par}
Note, from \rep{jacob} and \eqref{equil}, that when  $\frac{\delta^{1}}{\beta^{1}} = \frac{\delta^{2}}{\beta^{2}}$,
$$J(\tilde x^{1},\tilde x^{2})
\begin{bmatrix}
\phantom{-}\tilde x^{1}\; \cr
-\tilde x^{1}\; \end{bmatrix}=
0 \cdot
\begin{bmatrix}
\phantom{-}\tilde x^{1} \;\cr
-\tilde x^{1} \;\end{bmatrix},$$
i.e., the Jacobian matrix has a zero eigenvalue. Therefore, linearization says nothing about the local stability of the coexisting equilibria. Simulations indicate that, 
depending on the initial condition, the system can arrive at different equilibria of the form $\tilde x^{1} = \alpha\tilde x^{2}$ for different constants $\alpha>0$.
\hfill $\Box$
\end{remark}

Two viruses spreading on the same graph can be thought of as two products spreading in a market or two competing ideas spreading on a social network. Developing an understanding of how the viruses can coexist, and that the equilibrium that is reached is dependent on the initial condition, is vital to deploying initial marketing strategies that result in different market shares.  This coexistence  behavior is illustrated via simulation in Figures \ref{fig:par} and \ref{fig:pcomp}.

A similar result can be established for another special case, where two identical heterogeneous viruses spread on the same graph, as specified by the following assumption.

\begin{assumption}
Viruses 1 and 2 spread over the same strongly connected directed graph $\bbb{G}=([n],\scr{E})$, with  $\delta^{1}_i=\delta^{2}_i>0$ for all $i\in[n]$,
and $\beta^{1}_{ij}=\beta^{2}_{ij}$ for all $i=j\in[n]$ and $(i,j)\in\scr{E}$.
\label{special2}\end{assumption}

Under Assumption \ref{special2}, we have
$D^{1}=D^{2}=D$ and $B^{1}=B^{2}=B$,
where $D$ is a positive diagonal matrix and $B$ is an irreducible nonnegative matrix.

\begin{theorem}
Suppose that Assumptions \ref{x0} and \ref{special2} hold and that $s(-D+B)>0$.
If $(\tilde x^{1}, \tilde x^{2})$ with $\tilde x^{1},\tilde x^{2}>\0$ is an equilibrium
of \rep{sys}, then $\tilde x^{1},\tilde x^{2}\gg\0$, $\tilde x^{1}+\tilde x^{2}$ is unique, and
$\tilde x^{1} = \alpha\tilde x^{2}$ for some constant $\alpha>0$.
Furthermore, for each $\alpha>0$, there exists a unique pair $(\tilde x^{1}, \tilde x^{2})$ such that $\tilde x^{1} = \alpha\tilde x^{2}$.
\label{parallel2}\end{theorem}

{\em Proof:}
From \rep{sys} and Assumption \ref{special2},
\begin{eqnarray*}
\dot x^{1}(t) +\dot x^{2}(t) &=& \left(-D + B - (X^{1}(t)+X^{2}(t))B \right) \\
&& \ \ \ \ \ \ \ \times (x^{1}(t)+x^{2}(t)).
\end{eqnarray*}
Thus, the dynamics of $x^{1}(t)+x^{2}(t)$ is equivalent to that of the single-virus model \rep{single}.
By Proposition \ref{single0}, 
$x^{1}(t)+x^{2}(t)$ has a unique nonzero equilibrium in $[0,1]^n$.
Thus, $\tilde x^{1}+\tilde x^{2}$ is unique.
From \rep{sys}, we have
\begin{eqnarray*}
\dot x^{1}(t) -\dot x^{2}(t) &=& -D(x^{1}(t)-x^{2}(t)) +\\
&& (B - (X^{1}(t)+X^{2}(t))B) (x^{1}(t)-x^{2}(t)).
\end{eqnarray*}
Then,
$(-D+B - (\tilde X^{1}+\tilde X^{2})B) (\tilde x^{1}-\tilde x^{2}) =\0$.
Using the same arguments as those in the proof of Lemma \ref{nonzero},
it can be shown that $\tilde x^{1}+\tilde x^{2}\ll \1$. Thus,
$-D+B - (\tilde X^{1}+\tilde X^{2})B$ is an irreducible Metzler matrix.
By Lemma \ref{metzler}, since
$$(-D+B - (\tilde X^{1}+\tilde X^{2})B) (\tilde x^{1}+\tilde x^{2}) =\0,$$
and $x^{1}(t)+x^{2}(t)\gg \0$,
$s(-D+B - (\tilde X^{1}+\tilde X^{2})B)=0$.
From Lemma \ref{metzler0},  either $x^{1}(t)=x^{2}(t)$ or
$x^{1}(t)-x^{2}(t) = \gamma (x^{1}(t)+x^{2}(t))$
for some constant $\gamma>0$.
In both cases, it must be true that
$\tilde x^{1} = \alpha\tilde x^{2}$ for some constant $\alpha>0$,
and thus $\tilde x^{1},\tilde x^{2}\gg\0$.

Since $\tilde x^{1}+\tilde x^{2}$ is unique and $\tilde x^{1} = \alpha\tilde x^{2}$ for some constant $\alpha>0$, then
$\tilde x^{1}+\tilde x^{2} = (1 + \alpha )\tilde x^{2}$
is constant. Therefore, for each $\alpha>0$ there exists a unique $\tilde x^{2}$. So for each $\alpha>0$ there exists a unique pair $(\tilde x^{1}, \tilde x^{2})$.
\hfill
$\qed$



In this section, we have explored the equilibria of the bi-virus model and their local and global stability. It has been shown that unless both $s(- D^{1}+B^{1})$ and $s(- D^{2}+B^{2})$ are less than or equal to zero, at least one virus will pervade the network. Indeed, under different appropriate assumptions, one virus will fully dominate the competition, driving the second virus out, or both  viruses will coexist with infinitely many possible equilibria. 
It is worth noting that the necessary and sufficient 
condition for the eradication of both viruses (i.e., $s(- D^{1}+B^{1})\leq 0$ and $s(- D^{2}+B^{2})\leq 0$) requires global information on the network, which is inefficient and sometimes impossible. 
Subsequently, we explore some simple, but local, control techniques, with the expectation of attenuating or eliminating epidemic spreading.
We begin with the case when only one virus pervades.

\section{Sensitivity}\label{sense}

In this section, we regard each healing or infection rate as a local variable that an individual can control. The aim of this section is to understand whether and how local adjustment of healing and infection rates will affect the whole network.

We have already shown that in the case when $s(- D^{1}+B^{1})>0$ and $s(- D^{2}+B^{2})\le 0$,
the system \rep{sys} has a unique epidemic state of the form $(\tilde x^{1},\0)$ with $\tilde x^{1}\gg\0$, which is asymptotically stable.
It can be seen that the value of $\tilde x^{1}$ is independent of the matrices $D^{2}$ and $B^{2}$, but
depends on the matrices $D^{1}$ and $B^{1}$,
or equivalently, the parameters $\delta^{1}_i$ and $\beta^{1}_{ij}$.
A natural question is: how does the equilibrium $\tilde x^{1}$ change when the values of
$\delta^{1}_i$ and $\beta^{1}_{ij}$ are perturbed?
The aim of this section is to answer this question.

From the proof of Theorem \ref{eglobal},
the value of $\tilde x^{1}$ equals the unique epidemic state, denoted $x^*$, of the single-virus model \rep{single}
when $s(-D+B)>0$.
Thus, to answer the question just raised, it is equivalent to
study how the equilibrium $x^*$ changes when the values of
$\delta_i$ and $\beta_{ij}$ are perturbed.

For our purposes, we assume in this section that $\delta_i>0$ for all $i\in[n]$.
Then, by Proposition \ref{ff}, $s(-D+B)>0$ if and only if $\rho(D^{-1}B)>1$.

Suppose that $\rho(D^{-1}B)>1$. By Proposition \ref{single0}, the epidemic state $x^*$ is
the unique nonzero equilibrium of \rep{single}, which satisfies the equation
$(-D+B-X^*B)x^*=\0$.
Define the mapping $\Phi$ as follows:
$$\Phi(x^*,D,B):=(-D+B-X^*B)x^*.$$
Then, the equation $\Phi(x^*,D,B)= \0$ defines an implicit function $g:\R^{n\times n}\times \R^{n\times n}\rightarrow \R^n$
given by $x^*=g(D,B)$.

For each pair of matrices $D$ and $B$ for which $\rho(D^{-1}B)>1$,
there must exist a small neighborhood $\scr{B}$ such that for any pair of
matrices $D+\Delta D$ and $B+\Delta B$ in $\scr{B}$,
$$\rho\left((D+\Delta D)^{-1}(B+\Delta B)\right)>1.$$
Here $\Delta D$ is the $n\times n$ diagonal matrix whose $i$th diagonal entry equals $\Delta \delta_i$, which denotes the perturbation
of $\delta_i$, and
$\Delta B$ is the $n\times n$ matrix whose $ij$th entry equals $\Delta \beta_{ij}$, which denotes the perturbation of $\beta_{ij}$.
Let $x^*+\Delta x^*$ denote the new epidemic state resulting from the perturbations. Then,

\small \begin{equation*}
\left(-D-\Delta D+B+\Delta B - (X^*+\Delta X^*)(B+\Delta B)\right)
 (x^*+\Delta x^*) =\0,
\end{equation*}\normalsize
where $\Delta X^*={\rm diag}(\Delta x^*)$.
By ignoring the higher-order $\Delta$ terms, it is straightforward to verify that
\begin{equation}
\left(-D+B-X^*B-\tilde{B}^*\right)\Delta x^* \approx X^*\Delta \delta + (X^*-I)\Delta Bx^*,\label{dd}
\end{equation}
where $\Delta \delta$ is the vector in $\R^n$ whose $i$th entry equals  $\Delta \delta_i$ and $\tilde{B}^* = {\rm diag}(Bx^*)$.
Since we are interested in the local behavior around the equilibrium $x^*$, which is equivalent to performing linearization at $x^*$, the approximation in \rep{dd} is accurate for the following arguments.
First note that
$$\left(-D+B-X^*B-\tilde{B}^*\right)x^* = -\tilde{B}^*x^*.$$
Since $B$ is an irreducible nonnegative matrix and $x^*\gg \0$,
it follows that $\tilde{B}^*$ is a positive diagonal matrix.
Let $c>0$ be any positive constant, strictly smaller than
the minimal diagonal entry of $\tilde{B}^*$. Then,
$\tilde{B}^*> cI$ and thus $-\tilde{B}^*x^* < -cx^*$.
Since $(-D+B-X^*B-\tilde{B}^*)$ is an irreducible Metzler matrix,
by Lemma \ref{metzler},  $s(-D+B-X^*B-\tilde{B}^*)<-c<0$,
which implies that $(-D+B-X^*B-\tilde{B}^*)$ is nonsingular.
Thus, by the Implicit Function Theorem (see, e.g., pages 204-206 in \cite{chiang}),
 the function $x^*=g(D,B)$ is differentiable in the neighborhood $\scr{B}$.
From \rep{dd}, we have
\begin{eqnarray*}
\Delta x^* = \left(-D+B-X^*B-\tilde{B}^*\right)^{-1}X^*\Delta \delta +\;\;\;\;\;\;\;\;\;\;\;\;\\
\left(-D+B-X^*B-\tilde{B}^*\right)^{-1}(X^*-I)\Delta Bx^*.
\end{eqnarray*}

To proceed, we need the following lemma.

\begin{lemma}
{\rm (Theorem 2.7 in Chapter 6 of \cite{siam})}
Suppose that $M$ is a nonsingular, irreducible Hurwitz Metzler matrix.
Then, $M^{-1}\ll 0$.
\label{inverse}
\end{lemma}

From this lemma and the preceding discussion, it follows immediately that $(-D+B-X^*B-\tilde{B}^*)^{-1}$
is a strictly negative matrix.
Since $\0\ll x^* \ll \1$, it follows that
all $x_i^*$'s strictly decrease when any $\delta_i$ increases or any $\beta_{ij}$ decreases.
We have thus proved the following result.

\begin{proposition}
Consider the single-virus model \rep{single}.
Suppose that $\delta_i>0$ for all $i\in[n]$, and that the matrix $B$ is nonnegative and irreducible.
If $s(-D+B)>0$,
then each entry of the epidemic state $x^*$ is a strictly decreasing function of $\delta_i$, $i\in[n]$,
and a strictly increasing function of $\beta_{ij}$, $i,j\in[n]$.
\label{change1}
\end{proposition}

Similarly, we have the following result for the bi-virus model \rep{sys}.

\begin{theorem}
Suppose that $\delta^{1}_i>0$, $\delta^{2}_i\ge 0$ for all $i\in[n]$, and that matrices $B^{1}$ and $B^{2}$ are nonnegative and irreducible.
If $s(-D^{1}+B^{1})>0$ and $s(-D^{2}+B^{2})\le 0$,
then each entry of the epidemic state $\tilde x^{1}$ is a strictly decreasing function of $\delta^{1}_i$, $i\in[n]$,
and a strictly increasing function of $\beta^{1}_{ij}$, $i,j\in[n]$.
\label{change}
\end{theorem}

Theorem~\ref{change} characterizes the effects of adjusting each individual's healing rate and infection rates, which can be regarded as the simplest local control technique. Such adjustments can be achieved, for example, by taking medicine (i.e., increasing the healing rate) or reducing contact with neighbors (i.e., decreasing the infection rates). The result of Theorem~\ref{change} shows that any individual's local adjustment can attenuate every individual's epidemic state and thus the pervasion of the dominant virus.

It can be seen that if the individuals have sufficiently large healing rates and small infection rates, both $s(-D^1+B^1)$ and $s(-D^2+B^2)$ will be less than zero and thus both viruses will be eradicated. But to decide whether the healing rates are large enough or not, as well as whether the infection rates are small enough or not, requires centralized computation. With this in mind, we are interested in exploring distributed control techniques to eliminate (not only attenuate) epidemic spreading. It turns out that this is a challenging problem as discussed in the next section.

\section{Distributed Feedback Control}\label{feedback}

In this section, we regard each healing rate as a local control input of each group (or agent).
We begin with the single-virus model \rep{single}.

Suppose that the matrix $B$ is fixed. Let $\delta_i=\sum_{j=1}^n \beta_{ij}$.
Then, the row sums of $D^{-1}B$ all equal $1$.
By Lemma \ref{metzler0}, 
$\rho(D^{-1}B)=1$, which is equivalent to $s(-D+B)=0$ because of Proposition \ref{ff}. Thus, by Proposition \ref{single0}, the healthy state $\0$ is asymptotically stable
in this case.
This observation implies that in the case when local control inputs $\delta_i$'s are constant,
there always exist sufficiently large $\delta_i$'s which can stabilize the heathy state.

In the following, we will consider local control inputs of the form
\eq{\delta_i(t)=k_ix_i(t), \;\;\;\;\; i\in[n],\label{feed}}
where $k_i$ is a feedback gain. Designing the controller as a (linear) function of the infection proportion $x_i(t)$ is an intuitive approach since if the virus is eradicated, no control should be necessary. In implementation, this can be regarded as a treatment plan for individuals via administration of antidote or alternate treatment techniques.
By~\rep{update}, the system reduces to
\begin{eqnarray*}
\dot x_i(t) &=& -k_i (x_i(t))^2+(1-x_i(t))\sum_{j=1}^n \beta_{ij}x_j(t), \\
&& x_i(0)\in[0,1],\;\;\;\;\; i\in[n].\label{update3}
\end{eqnarray*}
The resulting $n$ state equations can be combined to yield 
\eq{\dot x(t) = \left(-KX(t)+B-X(t)B \right)x(t), \label{sys2}}
where $K= {\rm diag}([k_1,\dots , k_n])$.
Similar to the original system \rep{single}, we impose the following restrictions on the parameters of the new system \rep{sys2}.

\begin{assumption}
For all $i\in[n]$, we have $k_i>0$ and the matrix $B$ is nonnegative and irreducible.
\label{para2}
\end{assumption}


Using the same arguments as in the proof of Lemma \ref{box}, it is straightforward to
verify that the set $[0,1]^n$ is still positively invariant for the new system \rep{sys2}.
Since both $K$ and $X(t)$ are diagonal matrices, they commute. Then, from \rep{sys2},
\begin{eqnarray*}
\dot x(t) &=& (-X(t)K+B-X(t)B )x(t) \\
&=& (B-X(t)(K+B) )x(t) \\
&=& (-K+(K+B)-X(t)(K+B) )x(t).
\end{eqnarray*}
Thus, the system \rep{sys2} has the same form as the original system \rep{single},
with $D$ and $B$ being replaced by $K$ and $K+B$, respectively.

Note that $K^{-1}(K+B)=I+K^{-1}B$. Since by Assumption~\ref{para2}, $K$ is a positive diagonal matrix
and $B$ is an irreducible nonnegative matrix, $K^{-1}$ is a positive diagonal matrix and, thus,
$K^{-1}B$ is an irreducible nonnegative matrix.
By Lemma \ref{metzler0}, 
$\rho(K^{-1}B)>0$ and, thus, $\rho(I+K^{-1}B)>1$.
This observation implies, by Proposition \ref{2local},
that the new system \rep{sys2} has a unique nonzero equilibrium $x^*$ which satisfies
$\0\ll x^* \ll \1$ and is asymptotically stable with  domain of attraction
$[0,1]^n\setminus \{\0\}$.
We are thus led to the following result.

\begin{proposition}
Let Assumption \ref{para2} hold, and let $x(0)> \0$.
Then, for any local control inputs of the form
\rep{feed}, the healthy state $\0$ is not a reachable state of the system \rep{sys2}.
\label{im}
\end{proposition}

Note that instability of the healthy state can also be shown using Jacobian linearization at $\0$. However the above shows that the origin is not only an unstable equilibrium but actually is a repeller, that is, a perturbation in any direction will drive the state away from $\0$ toward equilibrium $x^*\gg \0$.

Now we turn to the bi-virus model \rep{sys}.
We consider local control inputs of the form
\eq{\delta^{1}_i(t)=k^{1}_ix^{1}_i(t), \;\;\;\;\; \delta^{2}_i(t)=k^{2}_ix^{2}_i(t), \;\;\;\;\; i\in[n],\label{feed2}}
where $k^{1}_i$ and $k^{2}_i$ are feedback gains.
By~\rep{updates}, the system reduces to
\begin{eqnarray*}
\dot{x}^{1}_i(t) &=& - k^{1}_i (x^{1}_i(t))^2 + (1 - x^{1}_i(t) - x^{2}_i(t))\sum_{j=1}^n \beta^{1}_{ij}  x^{1}_j(t) , \\
\dot{x}^{2}_i(t) &=& - k^{2}_i (x^{2}_i(t))^2 +(1 - x^{2}_i(t) - x^{1}_i(t))\sum_{j=1}^n \beta^{2}_{ij}  x^{2}_j(t) .
\end{eqnarray*}
The above equations can be written in matrix form:
\begin{equation}\label{sysfb}
\begin{split}
\dot{x}^{1}(t) &= (- K^{1}+(K^{1}+B^{1})-X^{1}(t)(K^{1}+B^{1}))x^{1}(t) \\
& \ \ \ \ \ \ \ \ \ \ \ \ \ \ \ \ \ \  \ \ \ \ \ \ \ \ \  - X^{2}(t)B^{1}x^{1}(t), \\
\dot{x}^{2}(t) &= (- K^{2}+(K^{2}+B^{2})-X^{2}(t)(K^{2}+B^{2}))x^{2}(t) \\
& \ \ \ \ \ \ \ \ \ \ \ \ \ \ \ \ \ \ \ \ \ \ \ \ \ \ \   - X^{1}(t)B^{2}x^{2}(t),
\end{split}
\end{equation}
where $K^{1}$ and $K^{2}$ are $n\times n$ diagonal matrices with the $i$th diagonal entries equal to $k^{1}_i$
and $k^{2}_i$, respectively.
Similar to the original system \rep{sys}, we impose the following restrictions on the parameters of the new system \rep{sysfb}.

\begin{assumption}
For all $i\in[n]$, we have $k^{1}_i,k^{2}_i>0$ and the matrices $B^{1}$ and $B^{2}$ are nonnegative and irreducible.
\label{para3}
\end{assumption}

From the preceding discussion, we have the following.

\begin{theorem}
Let Assumptions \ref{x0} and \ref{para3} hold.
Then, for any local control inputs of the form
\rep{feed2}, the healthy state $(\0,\0)$ is an unstable equilibrium of the system \rep{sysfb}.
\label{im2}
\end{theorem}

Theorem~\ref{im2} implies that the distributed feedback controller \rep{feed2} can never stabilize the healthy state.  
This impossibility result is interesting and surprising because proportional controllers, while suboptimal, function fairly well in many applications. The result, as well as Proposition~\ref{im}, thus partially explains why distributed control of epidemic networks is a challenging open problem.


\section{Simulations}\label{simulations}

In this section we first compare the model in \eqref{sys} to the full probabilistic $3^n$ state model in  \eqref{eq:qij}-\eqref{eq:v}. We then illustrate some of the results from Section \ref{equilibria} via simulation. Finally we explore some interesting behavior via simulations. 

\begin{figure}
\centering
\includegraphics[width=.8\columnwidth]{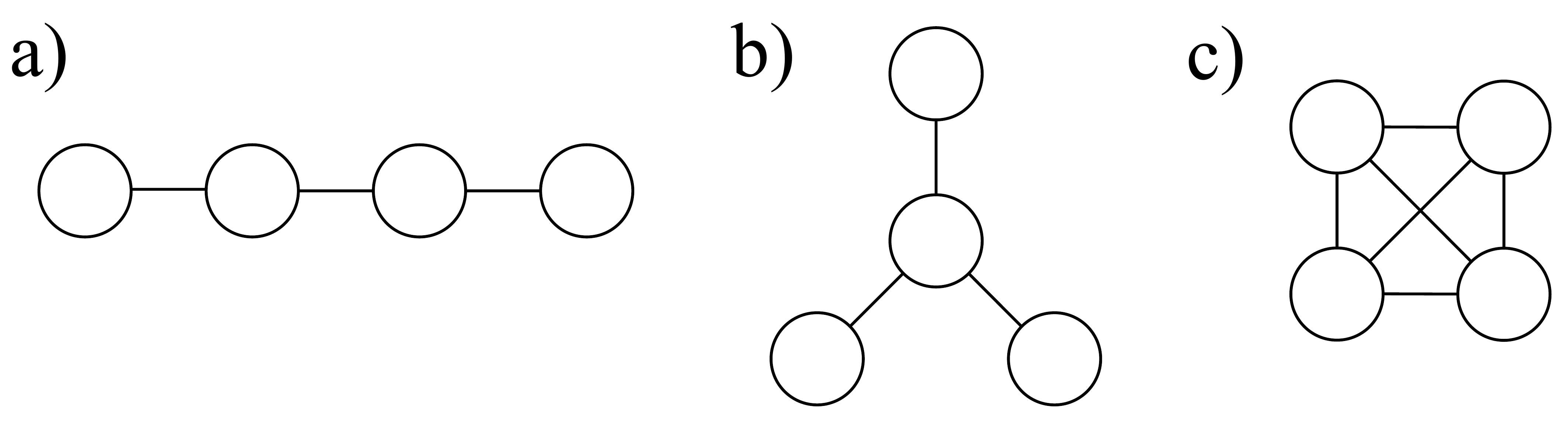}
\caption{Graph structures: a) line, b) star, c) complete.}
\label{fig:graphs}
\end{figure}

\begin{figure}
    \centering
\begin{overpic}[width=\columnwidth]{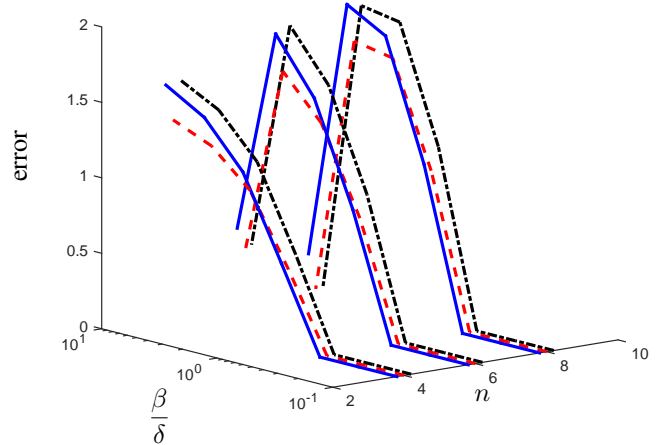}
    \put(69,6){{\parbox{\linewidth}{%
     $n$}}}
     \put(20,3){{\parbox{\linewidth}{%
     $\displaystyle\frac{\beta}{\delta}$}}}
     \put(0,42){{\parbox{\linewidth}{%
     \rotatebox{90}{error}}}}
\end{overpic}
\caption{
A plot of $\|[v^1(T); v^2(T)] -[x^1(T); x^2(T)]\|$ for the line graph, $T=10,000$. Results from using the different initial conditions $(x^{1_1}(0),x^{2_1}(0))$, $(x^{1_2}(0),x^{2_2}(0))$, $(x^{1_3}(0),x^{2_3}(0))$ are depicted by the blue lines, red dashed lines, and black dash-dot lines, respectively. 
}
\label{fig:line}
\end{figure}

\subsection{Approximation Accuracy}

We evaluate the effectiveness of  \eqref{sys} (with $\beta_{ii}=0 \ \forall i$) as an approximation of the $3^n$ state model in \eqref{eq:qij}-\eqref{eq:v} for line graphs, star (hub--spoke) graphs, and complete graphs; see Figure \ref{fig:graphs} for examples of each graph structure. 
All adjacency matrices for these graphs are symmetric and binary-valued, and both viruses spread on the same graph. In the star graph, the central node is the first agent. 
Each simulation was run for 10,000 time steps  (final time $T=10,000$), with three initial conditions: 1) the first node infected with virus 1 and the second node infected with virus 2, $x^{1_1}(0) = [ 1 \  0 \cdots\ 0]^{\top}, x^{2_1}(0) = [ 0 \ 1 \ 0 \ \cdots\ 0]^{\top}$, 2) the first two nodes infected with virus 1 and the second two nodes infected with virus 2, $x^{1_2}(0) = [ 1 \ 1 \ 0 \cdots\ 0]^{\top}, x^{2_2}(0) = [ 0 \ 0 \ 1 \ 1 \ 0 \ \cdots\ 0]^{\top}$, and 3) the first node infected with virus 1 and the rest of the  nodes infected with virus 2, $x^{1_3}(0) = [ 1  \ 0 \cdots\ 0]^{\top}, x^{2_3}(0) = [ 0  \ 1 \ \cdots\ 1]$. 
We  explore identical homogeneous viruses, $(\beta,\delta) = (\beta^1,\delta^1) = (\beta^2,\delta^2)$ in these tests. The $(\beta,\delta)$ pairs used are $[ (0.1,1), (0.215,1), (0.464,1),(0.5,0.5), (1,0.464),$ $ (1,0.215), (1, 0.1)]$, and the number of agents, $n = 4,6,8$. 
We limited simulations to these $n$ values since mean field approximations are typically worse for small values of $n$ and there is a computational limitation due to the size of the $3^n$ state model.

\begin{figure}
    \centering
\begin{overpic}[width=\columnwidth]{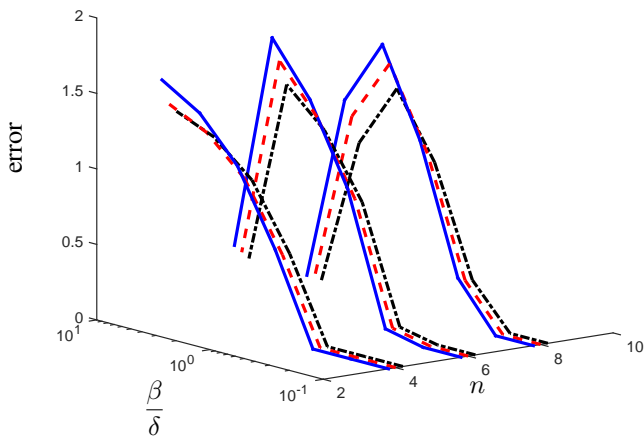}
    \put(69,6){{\parbox{\linewidth}{%
     $n$}}}
     \put(20,3){{\parbox{\linewidth}{%
     $\displaystyle\frac{\beta}{\delta}$}}}
     \put(0,42){{\parbox{\linewidth}{%
     \rotatebox{90}{error}}}}
\end{overpic}
\caption{
A plot of $\|[v^1(T); v^2(T)] -[x^1(T); x^2(T)]\|$ for the star graph, $T=10,000$. Results from using the different initial conditions $(x^{1_1}(0),x^{2_1}(0))$, $(x^{1_2}(0),x^{2_2}(0))$, $(x^{1_3}(0),x^{2_3}(0))$ are depicted by the blue lines, red dashed lines, and black dash-dot lines, respectively. 
}
\label{fig:star}
\end{figure}

\begin{figure}
    \centering
    \begin{overpic}[width=\columnwidth]{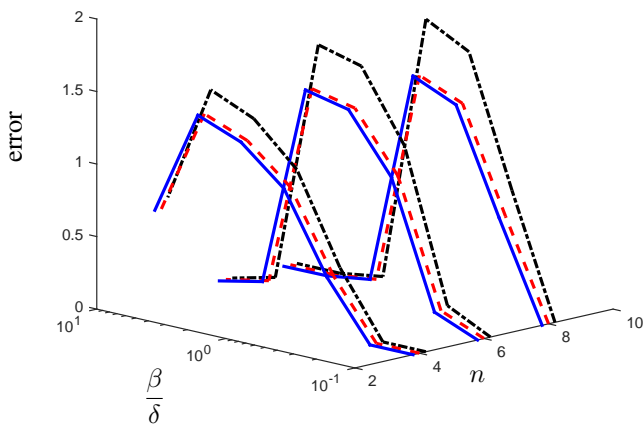}
    \put(69,6){{\parbox{\linewidth}{%
     $n$}}}
     \put(20,3){{\parbox{\linewidth}{%
     $\displaystyle\frac{\beta}{\delta}$}}}
     \put(0,42){{\parbox{\linewidth}{%
     \rotatebox{90}{error}}}}
\end{overpic}
\caption{
A plot of the error $\|[v^1(T); v^2(T)] -[x^1(T); x^2(T)]\|$ for the complete graph, $T=10,000$. Results from using the different initial conditions $(x^{1_1}(0),x^{2_1}(0))$, $(x^{1_2}(0),x^{2_2}(0))$, $(x^{1_3}(0),x^{2_3}(0))$ are depicted by the blue lines, red dashed lines, and black dash-dot lines, respectively. 
}
\label{fig:full}
\end{figure}

The results are given in Figures \ref{fig:line}-\ref{fig:full} in terms of the 2-norm of the difference between the states of \eqref{sys} at the final time ($[x^1(T); x^2(T)]$), and the means of the two viruses in the $3^n$ state Markov model at the final time ($[v^1(T); v^2(T)]$ as defined by \eqref{eq:v}). 

The accuracy of the approximation appears to be very similar to the single virus case \cite{OmicTN09,pare2017epidemic}. 
Since \eqref{sys} is an upper bounding approximation, the results show that the two models converge to the healthy state for the smaller values of $\frac{\beta }{\delta}$, resulting in small errors between the two models. For many of the larger values of $\frac{\beta }{\delta}$, \eqref{sys} again performs quite well since it is at an epidemic state and the $3^n$ state model does not appear to reach the healthy state in the finite time considered in the simulations ($T=10,000$). Therefore for certain values of $\frac{\beta }{\delta}$ and certain time scales, \eqref{sys} is a sufficient approximation of the $3^n$ state Markov model. 
For values of $\frac{\beta }{\delta}$ that are near one, the models are quite different, similar to the single virus case. The $3^n$ state model appears, in most cases, to be at or close to the healthy state while \eqref{sys} is at an epidemic state, resulting in large errors. 


\subsection{Illustrative Examples}


\begin{figure}
    \centering
    \begin{subfigure}[b]{.4935\columnwidth}
      \includegraphics[width=\columnwidth]{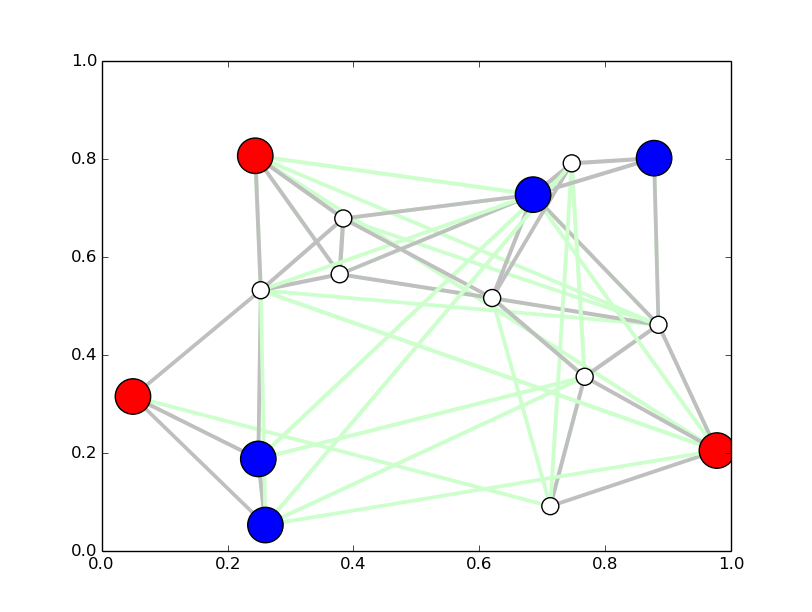}
      \caption{The system at time zero.}
      \label{fig:dfe0}
    \end{subfigure}
    \hfill
    \begin{subfigure}[b]{.4935\columnwidth}
      \includegraphics[width=\columnwidth]{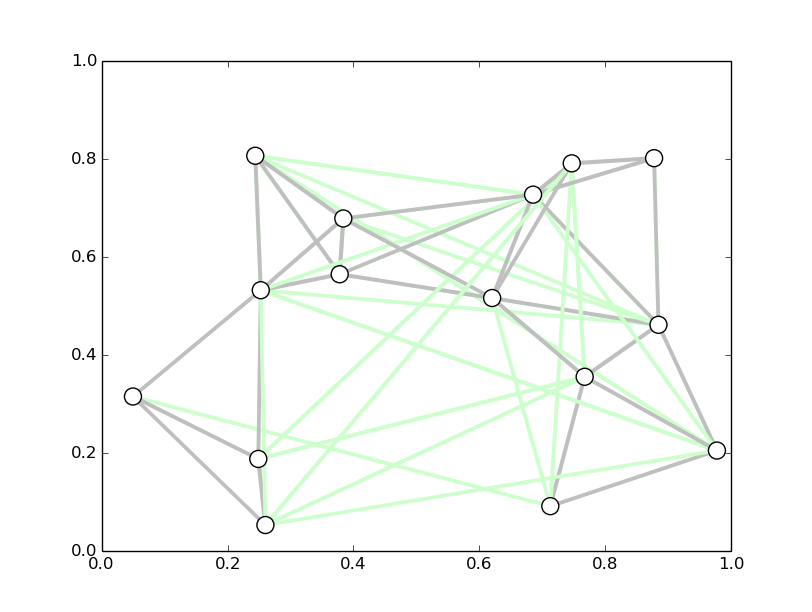}
      \caption{The system at time 400.}
      \label{fig:dfe400}
    \end{subfigure}
    \caption{This bi-virus system meets the assumptions of Theorem~\ref{0global}, with $s(-D^{1}+B^{1})=-0.1191$ and $s(-D^{2}+B^{2})=-0.0316$, and both viruses are eradicated in 400 time steps.  
    The colors and diameters follow \eqref{eq:color} and \eqref{eq:diam}.
For a video of this simulation please see \href{http://youtu.be/ZPjk52uiJi0}{youtu.be/ZPjk52uiJi0}.}
\label{fig:dfe}
\end{figure}

\begin{figure}
\centering
\includegraphics[width = .8\columnwidth]{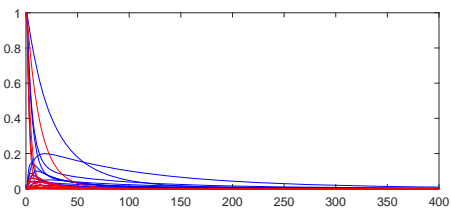}
\caption{
Trajectories of the simulation in Figure \ref{fig:dfe}. Virus 1 and 2 are depicted by red and blue, respectively.}
\end{figure}
This section contains several illustrative and insightful simulations. 
Given the nature of the network-dependent competing viruses and the behavior of the parallel equilibria ($\tilde x^{1} = \alpha\tilde x^{2}$) shown by Theorems \ref{parallel1} and \ref{parallel2}, we employ a unique coloring scheme. In each figure (Figures~4-6), we plot the initial condition on the left and the final state on the right, and we include a link to a video of the full simulation in the caption. The trajectories of the two viruses for each node are then depicted by the color of the corresponding node.
Virus 1 is depicted by the color red ($r= [ 1 \ 0 \ 0]$) and virus 2 is depicted by the color blue ($b= [ 0 \ 0 \ 1]$). For each $i\in[n]$, the color at time $t$ for group (agent) $i$ is 
\begin{equation}\label{eq:color}
\frac{x^{1}_i(t)}{x^{1}_i(t)+x^{2}_i(t)}r + \frac{x^{2}_i(t)}{x^{1}_i(t)+x^{2}_i(t)}b.
\end{equation}
When $x^{1}_i(t)+x^{2}_i(t)=0$, the color goes to white, indicating a completely healthy, susceptible group $i$.
These colors are used to facilitate the depiction of the parallel equilibria, which will be shown by all nodes converging to the same color.
For all $i\in[n]$, the diameter of the node representing group (agent) $i$ is given by
\begin{equation}\label{eq:diam}
d_0 + (x^{1}_i(t)+x^{2}_i(t))r_0,
\end{equation}
with $d_0$ being the default/smallest diameter and $r_0$ being the scaling factor depending on the total sickness of group (agent) $i$.
Therefore, the color indicates the proportion of each virus the group (agent) has and the diameter indicates the strength of the viruses or how sick the group (agent) is.

For systems with two different graph structures, the graph over which virus 1 spreads is depicted by gray edges and the graph over which virus 2 spreads is depicted by green edges.
If both viruses spread over the same graph, the edges are gray.

First we illustrate Theorem \ref{0global}, with $s(-D^{1}+B^{1})=-0.1191$ and $s(-D^{2}+B^{2})=-0.0316$.
See Figure \ref{fig:dfe} for the initial and final states.
Consistent with the result of the theorem, both viruses are eradicated.

\begin{figure}
    \centering
    \begin{subfigure}[b]{.493\columnwidth}
      \includegraphics[width=\columnwidth]{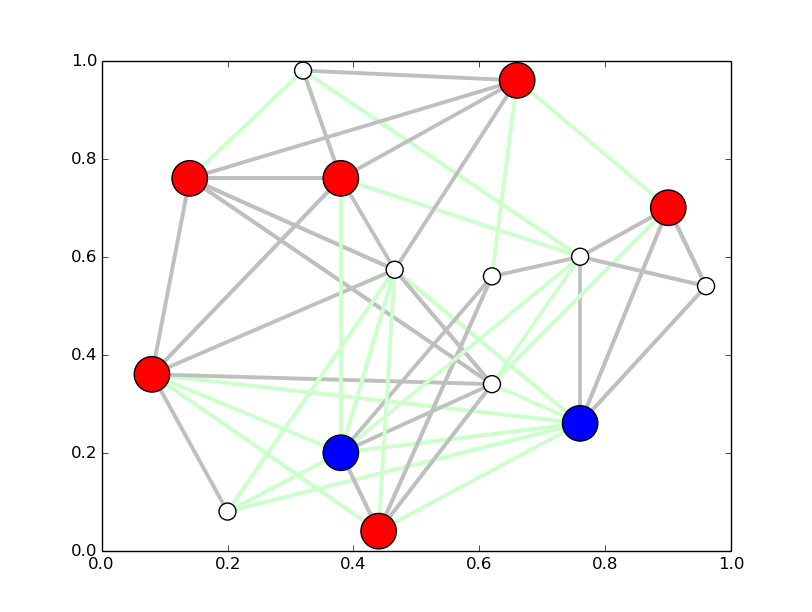}
      \caption{The system at time zero.}
      \label{fig:1dom0}
    \end{subfigure}
    \hfill
    \begin{subfigure}[b]{.493\columnwidth}
      \includegraphics[width=\columnwidth]{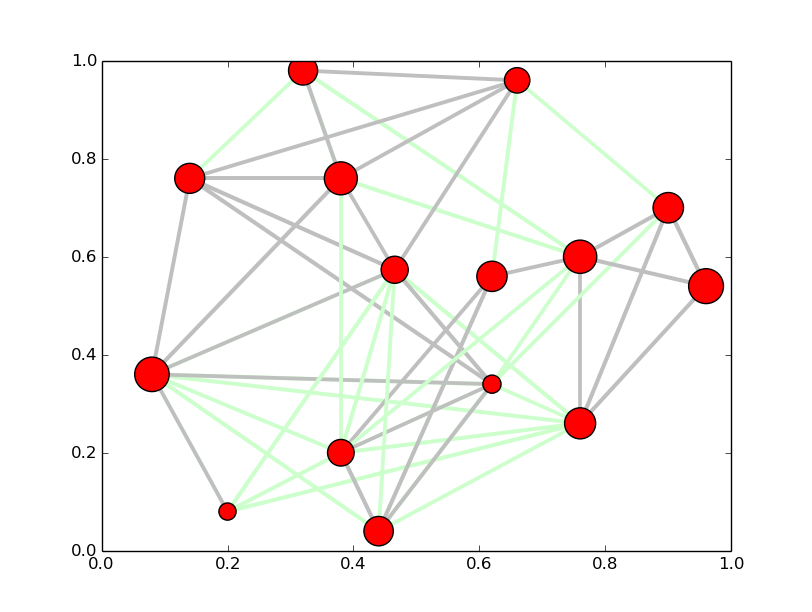}
      \caption{The system at time 300.}
      \label{fig:1dom400}
    \end{subfigure}
    \caption{This bi-virus system meets the assumptions of Theorem~\ref{eglobal}, with $s(-D^{1}+B^{1})=0.4145$ and $s(-D^{2}+B^{2})=-0.0802$. Virus 2 is eradicated and virus 1 reaches its epidemic state in 300 time steps.  
    The colors and diameters follow \eqref{eq:color} and \eqref{eq:diam}.
For a video of this simulation please see \href{http://youtu.be/dSm9P0O3c6A}{youtu.be/dSm9P0O3c6A}.}
\label{fig:1dom}
\end{figure}

We  illustrate Theorem \ref{eglobal}, with $s(-D^{1}+B^{1})=0.4145$ and $s(-D^{2}+B^{2})=-0.0802$,
with Figure \ref{fig:1dom} depicting the initial and final states.
Consistent with the result of the theorem, one virus reaches an epidemic state, while the other is eradicated.

\begin{figure}
\centering
\includegraphics[width = .8\columnwidth]{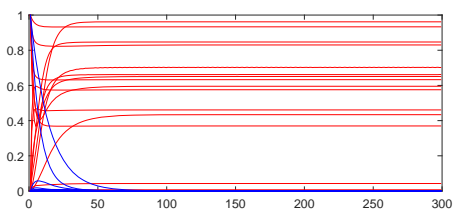}
\caption{
Trajectories of the simulation in Figure \ref{fig:1dom}. Virus 1 and 2 are depicted by red and blue, respectively.}
\end{figure}

\begin{figure}
    \centering
    \begin{subfigure}[b]{.493\columnwidth}
      \includegraphics[width=\columnwidth]{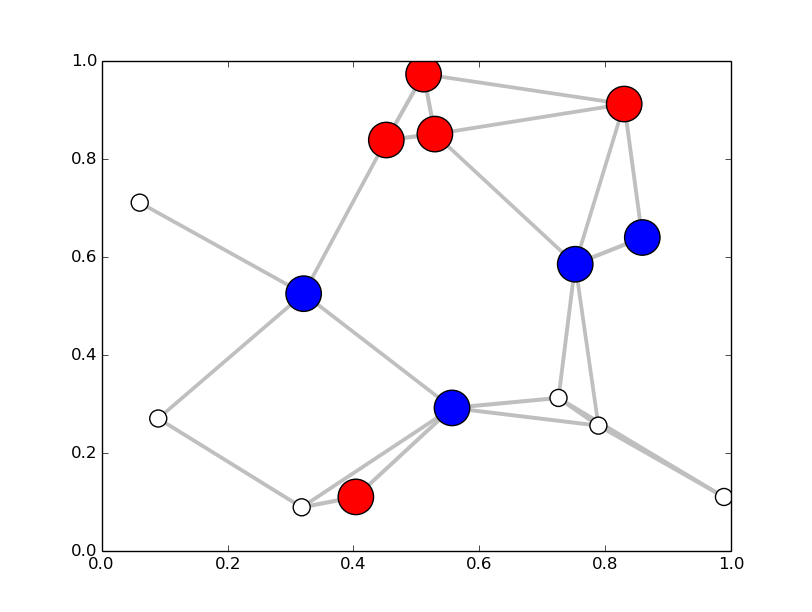}
      \caption{The system at time zero.}
      \label{fig:par0}
    \end{subfigure}
    \hfill
    \begin{subfigure}[b]{.493\columnwidth}
      \includegraphics[width=\columnwidth]{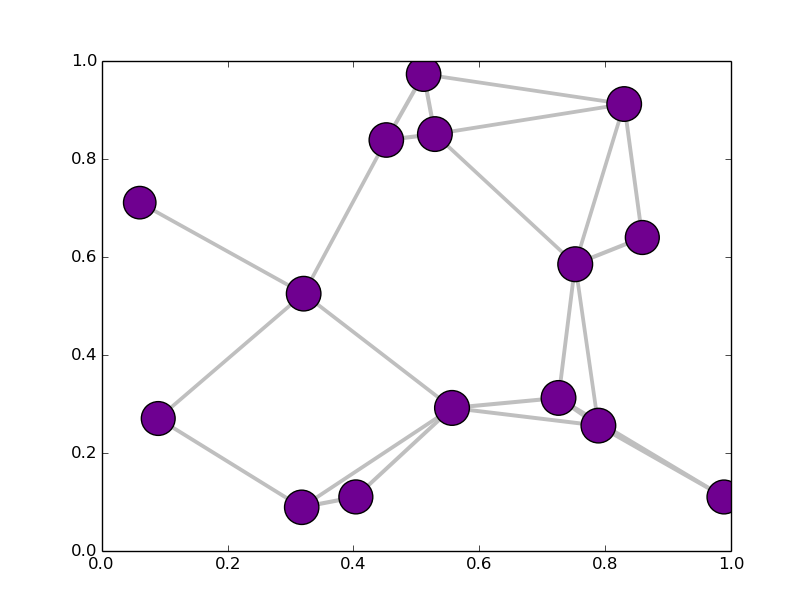}
      \caption{The system at time 1,000.}
      \label{fig:par400}
    \end{subfigure}
    \caption{This bi-virus system meets the assumptions of Theorem~\ref{homo1}, with $\frac{\delta^{1}}{\beta^{1}}=\frac{0.4242}{1.9090} = \frac{\delta^{2}}{\beta^{2}}=\frac{0.2121}{0.9545}$ and $s(-\delta^{1}I+\beta^{1})=4.2654 = 2s(-\delta^{2}I+\beta^{2}A)$, and the viruses converge to a parallel equilibrium.  
    The colors and diameters follow \eqref{eq:color} and \eqref{eq:diam}.
For a video of this simulation please see \href{http://youtu.be/Ui1BA3C0zI0}{youtu.be/Ui1BA3C0zI0}.}
\label{fig:par}
\end{figure}

Theorem \ref{parallel1} is illustrated in Figure \ref{fig:par}, with $\frac{\delta^{1}}{\beta^{1}}=\frac{0.4242}{1.9090} = \frac{\delta^{2}}{\beta^{2}}=\frac{0.2121}{0.9545}$ and $s(-\delta^{1}I+\beta^{1})=4.2654 = 2s(-\delta^{2}I+\beta^{2}A)$. 
Consistent with the result of the theorem, both viruses reach an epidemic state where one equilibrium is a scaling of the other, depicted by all the colors being the same, by \eqref{eq:color}.

\begin{figure}
    \centering
    \begin{subfigure}[b]{\columnwidth}
    \centering
      \includegraphics[width=.9\columnwidth]{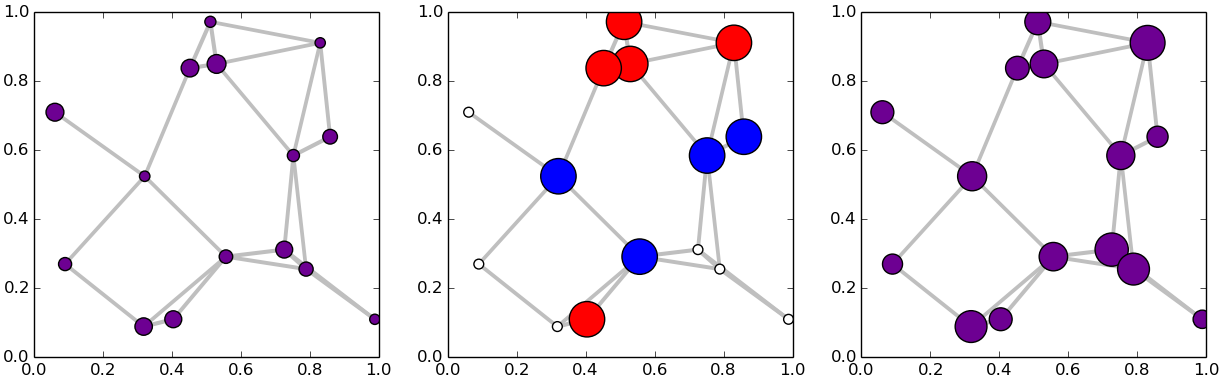}
      \caption{The system at time zero.}
      \label{fig:p0}
    \end{subfigure}
    \hfill
    \begin{subfigure}[b]{\columnwidth}
    \centering
      \includegraphics[width=.9\columnwidth]{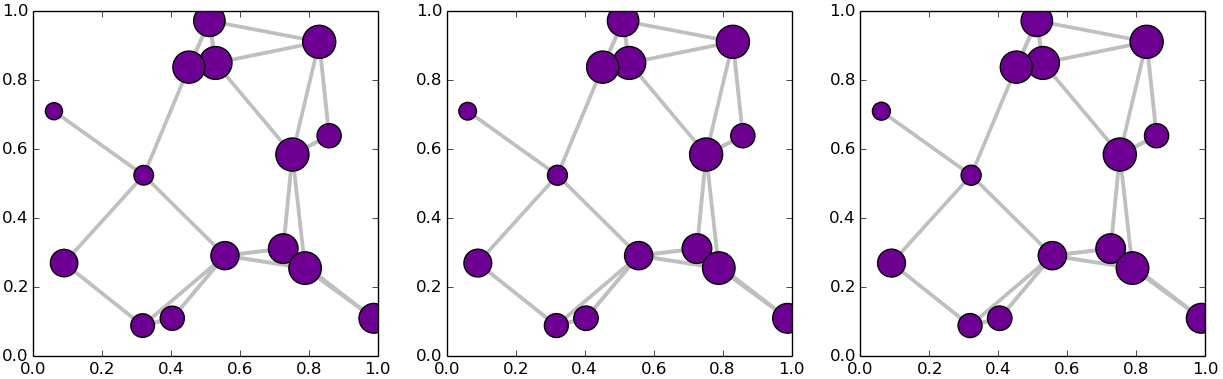}
      \caption{The system at time 1,000.}
      \label{fig:p400}
    \end{subfigure}
    \caption{This bi-virus system meets the assumptions of Theorem \ref{homo1} and the viruses converge to the same parallel equilibrium for all three initial conditions.  
    The colors and diameters follow \eqref{eq:color} and \eqref{eq:diam}.
For a video of this simulation please see \href{http://youtu.be/ES1MT-kWgIM}{youtu.be/ES1MT-kWgIM}.}
\label{fig:pcomp}
\end{figure}

A question of interest is, when do certain systems converge to different, or the same, parallel equilibrium? In the following simulation we start the system from Figure \ref{fig:par} with three different initial conditions and they all converge to the same equilibrium. This is depicted in Figure \ref{fig:pcomp}. 
However, consistent with Theorem \ref{parallel1}, the system has many equilibria which, via simulation, appear to be initial condition dependent.

\begin{figure}
\centering
\includegraphics[width = .32\columnwidth]{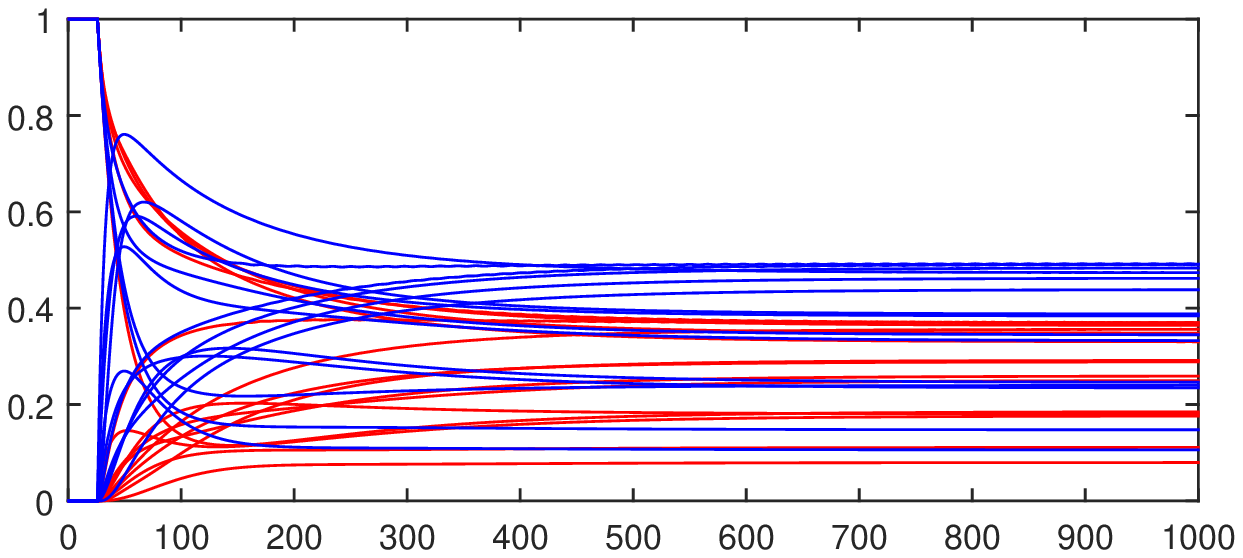} \includegraphics[width = .32\columnwidth]{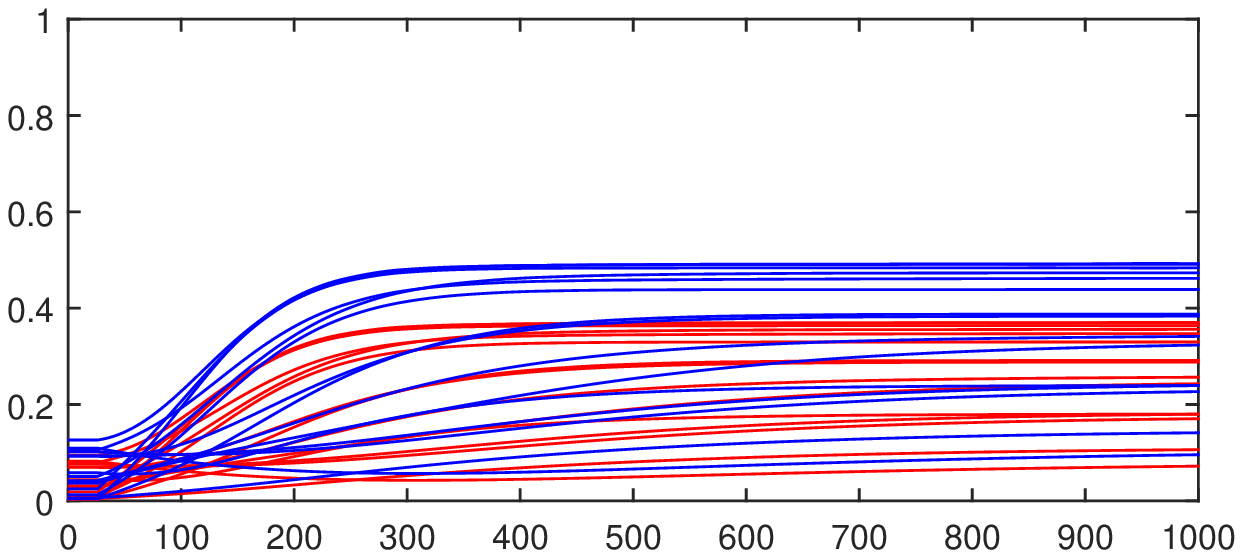} \includegraphics[width = .32\columnwidth]{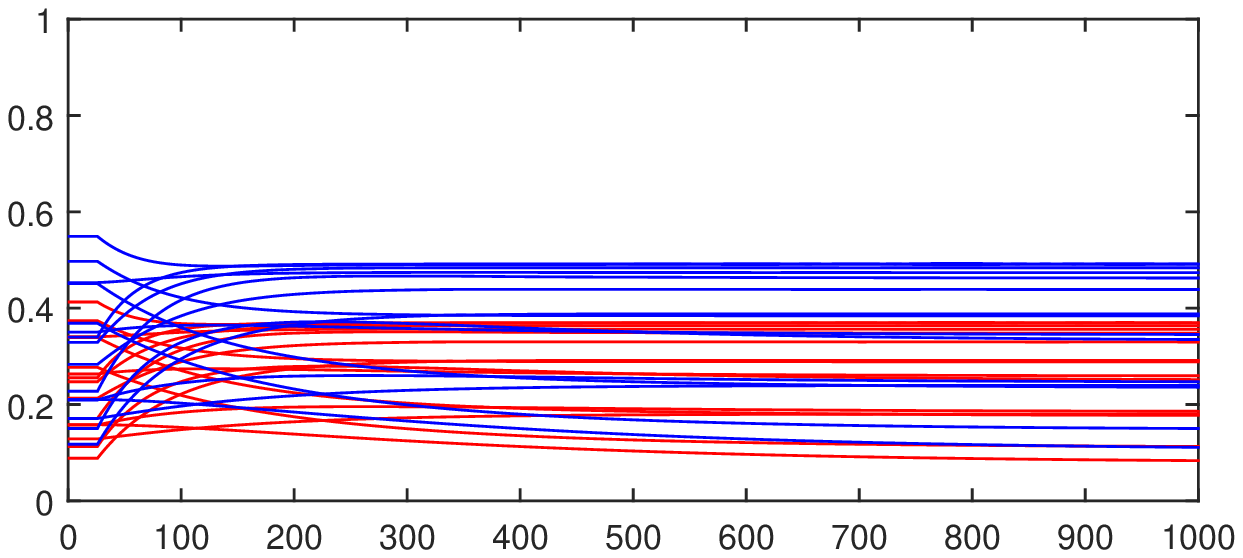}
\caption{
Trajectories of the simulations in Figure \ref{fig:pcomp}. Virus 1 and 2 are depicted by red and blue, respectively.}
\end{figure}

\begin{figure}
    \centering
    \begin{subfigure}[b]{.493\columnwidth}
      \includegraphics[width=\columnwidth]{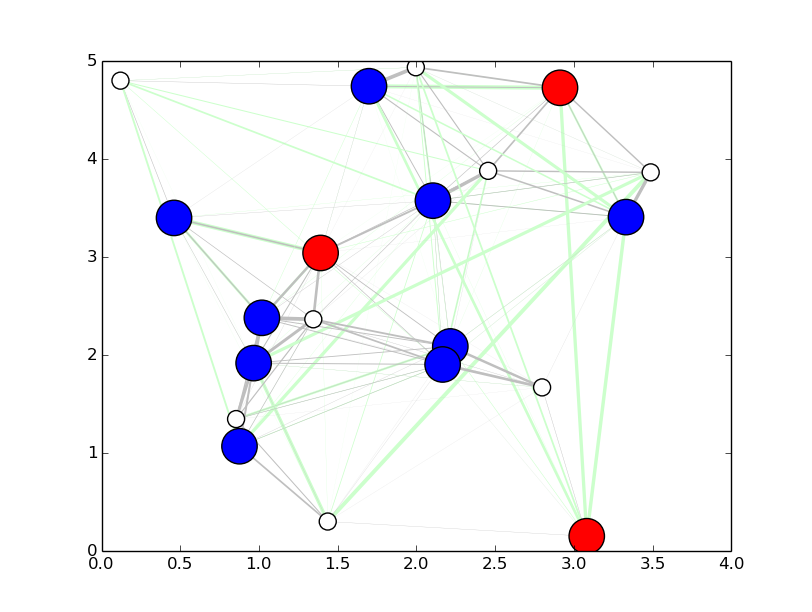}
      \caption{The system at time zero.}
      \label{fig:co0}
    \end{subfigure}
    \hfill
    \begin{subfigure}[b]{.493\columnwidth}
      \includegraphics[width=\columnwidth]{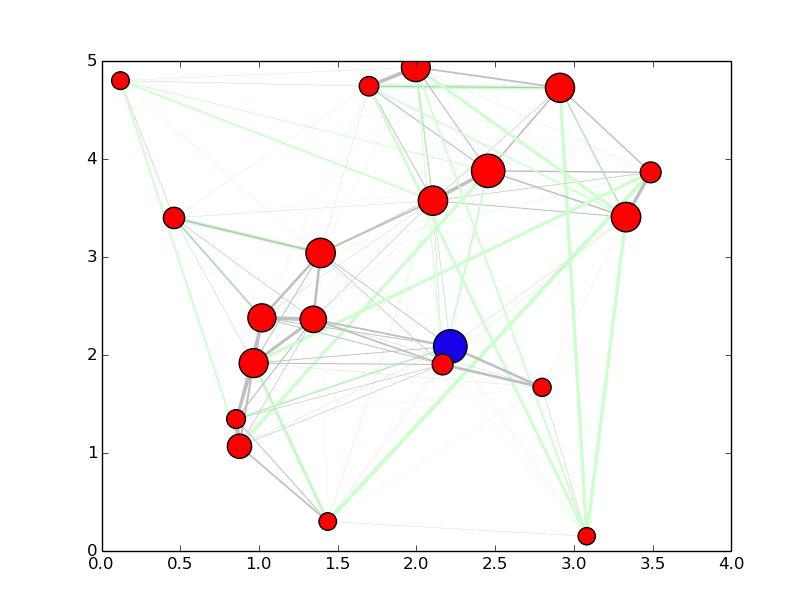}
      \caption{The system at time 400.}
      \label{fig:co400}
    \end{subfigure}
    \caption{This bi-virus system has $s(-D^{1}+B^{1})=0.2276> 0$, $s(-D^{2}+B^{2})=0.3117>0$, $D^{1}\neq D^{2}$, $B^{1} \neq B^{2}$, with virus 1 spreading on the gray edges, and virus 2 spreading on the green edges.  
    The colors and diameters follow \eqref{eq:color} and \eqref{eq:diam}.
For a video of this simulation please see \href{http://youtu.be/MTKoqposczA}{youtu.be/MTKoqposczA}.}
\label{fig:coexist}
\end{figure}

Finally, there was earlier discussion about coexisting equilibria and it was shown, for particular cases, that convergence to such equilibria can be proven. However simulations show that generic systems can also converge to epidemic equilibria. Consider the system in Figure \ref{fig:coexist}, where $s(-D^{1}+B^{1})=0.2276> 0$, $s(-D^{2}+B^{2})=0.3117>0$, $D^{1}\neq D^{2}$, $B^{1} \neq B^{2}$, and the viruses spread over different graphs.
Both viruses clearly reach an epidemic equilibrium; see Figure \ref{fig:co400}. Simulations show for the same system, that different initial conditions lead to an epidemic state for virus 1 and a healthy state for virus 2.
Therefore, in the generic case the system is initial condition dependent.

\begin{figure}
\centering
\includegraphics[width = .8\columnwidth]{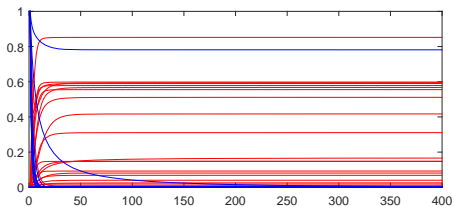}
\caption{
Trajectories of the simulation in Figure \ref{fig:coexist}. Virus 1 and 2 are depicted by red and blue, respectively.}
\end{figure}


\section{Conclusion}\label{conclusion}


In this paper, we have thoroughly analyzed the equilibria of a continuous-time bi-virus model and their stability, and in so doing, as a by-product we have  improved on the results for the single-virus model. We have provided necessary and sufficient conditions for convergence to the healthy state of \eqref{sys}. We have also provided  results on the epidemic states of \eqref{sys}, including several sufficient conditions for stability and instability.
We have explored two distributed control techniques for the model. First, we have regarded each healing or infection rate as a local variable that an individual can control, and shown by sensitivity analysis that any small local adjustment of healing or infection rates can attenuate the global pervasion of the dominant virus. Second, we have regarded each healing rate as a local control input, and 
shown that a  distributed proportional controller of the form $\delta_i(t)=k_ix_i(t)$, can never drive the virus model to the healthy state.

For future work, we plan to study bi-virus models with time-varying graph structures, similar to \cite{pare2015stability,pare2017epidemic} for the single-virus case. We also plan to analyze the multi-virus case, i.e., more than two competing viruses, for the healthy and epidemic states. 
The sensitivity analysis over the
regime $\rho(D_1^{-1}B_1)>0$ and $\rho(D_2^{-1}B_2) > 0$, where both strains would be fit to survive, constitutes a future research direction. 
Distributed control of epidemic networks, including the bi-virus model considered here, is of course another direction for future research.

We note that there are two limitations of the bi-virus model considered in the paper. First, it was assumed that individuals are susceptible at
birth, even if their parents are infected with a virus. This assumption is invalid
in some realistic situations as individuals can be infected with a virus at birth if their parents are
virus carriers. Second, the infection rate within the $i$th group ($\beta_{ii}$)
is probably larger than the infection rate between two different groups ($\beta_{ij} \ i \neq j$) in some realistic situations. 
Modeling and investigating the dynamics of bi-virus spreading in such scenarios are left as an important and interesting area for future research.


\section{Appendix}


{\em Proof of Lemma~\ref{pattern}:}
Suppose, to the contrary, that for all $i\in[n]$ such that $x_i=0$, we have $y_i=0$.
Since $x$ has at least one zero entry and $x>\0$, there exists a proper nonempty subset
$\scr{E}\subset [n]$ such that for  $i\in\scr{E}$, $x_i=0$,
and for $i\in [n] \setminus \scr{E}$, $x_i>0$. 
By our assumption, it follows that for any $i\in\scr{E}$, $y_i=0$.
Without loss of generality,
let $\scr{E}=\{1,2,\ldots,m\}$ for some $m$,  where $1\le m<n$.
Then,
$$Mx = \begin{bmatrix}
A & B \cr
C & D\end{bmatrix}
\begin{bmatrix}\; \0 \; \cr \bar x\end{bmatrix}
= y = \begin{bmatrix} \; \0 \; \cr \bar y\end{bmatrix},$$
with $A,D\ge0$, $B,C>0$, $\bar x \gg \0$, and $\bar y \ge \0$.
This implies the $B = 0$, which is a contradiction since $M$ is an irreducible matrix.
Therefore, there exists an index $i\in[n]$ such that $x_i=0$ and $y_i>0$,
and thus $x$ and $y$ cannot have the same sign pattern.
\hfill
$\qed$

{\em Proof of Proposition~\ref{ff}:}
Suppose that $\Lambda$ is a negative diagonal matrix in $\R^{n\times n}$
and $N$ is an irreducible nonnegative matrix in $\R^{n\times n}$.
Let $M=\Lambda+N$. By Theorem 3.29 in \cite{varga},
$s(M)<0$ if and only if $\rho(-\Lambda^{-1}N)<1$.
To prove the proposition, it suffices to show that $s(M)=0$ if and only if $\rho(-\Lambda^{-1}N)=1$.

First suppose that $s(M)=0$.
Set $\Lambda_{\varepsilon}=\Lambda - \varepsilon I$ with $\varepsilon>0$.
Let
$M_{\varepsilon} = \Lambda_{\varepsilon} +N = \Lambda - \varepsilon I +N$.
Then, $\lim_{\varepsilon\rightarrow 0^+}\rho(-\Lambda_{\varepsilon}^{-1}N) = \rho(-\Lambda^{-1}N)$.
Since $\varepsilon>0$, $s(M_{\varepsilon})<0$. Then,
$\rho(-\Lambda_{\varepsilon}^{-1}N) < 1$ and, therefore, $\lim_{\varepsilon\rightarrow 0^+}\rho(-\Lambda_{\varepsilon}^{-1}N) \le 1$.
Thus, $\rho(-\Lambda^{-1}N)\le 1$. To prove that $\rho(-\Lambda^{-1}N) = 1$,
suppose that, to the contrary, $\rho(-\Lambda^{-1}N)<1$.
Then, $s(M)<0$, which is a contradiction. Therefore, $\rho(-\Lambda^{-1}N) = 1$.

Now suppose that $\rho(-\Lambda^{-1}N) = 1$. Again set $\Lambda_{\varepsilon}=\Lambda - \varepsilon I$ with $\varepsilon>0$
and $M_{\varepsilon} = \Lambda_{\varepsilon} +N$.
Then, $\lim_{\varepsilon\rightarrow 0^+} s(M_{\varepsilon}) = s(M)$.
Since $\varepsilon>0$, $-\Lambda_{\varepsilon}^{-1}N$ is a nonnegative matrix.
Since $N$ is irreducible and nonnegative, so is $-\Lambda_{\varepsilon}^{-1}N$.
Note that the $i$th diagonal entry of $-\Lambda_{\varepsilon}$ is strictly larger than
the $i$th diagonal entry of $-\Lambda$ since $\varepsilon>0$. Thus,
$-\Lambda^{-1}N > -\Lambda_{\varepsilon}^{-1}N$.
By Lemma \ref{metzler0}, 
$\rho(-\Lambda_{\varepsilon}^{-1}N) < 1$.
Then,
$s(M_{\varepsilon})<0$ and, thus, $\lim_{\varepsilon\rightarrow 0^+} s(M_{\varepsilon})\le 0$.
Thus, $s(M)\le 0$. To prove that $s(M) = 0$,
suppose that, to the contrary, $s(M)<0$.
Then, $\rho(-\Lambda^{-1}N)<1$, which is a contradiction. Therefore, $s(M)=0$.
\hfill
$\qed$

{\em Proof of Lemma~\ref{box}:}
Suppose that at some time $\tau$, $x^{1}_i(\tau),x^{2}_i(\tau),x^{1}_i(\tau)+x^{2}_i(\tau)\in[0,1]$ for all $i\in[n]$.
Consider an index $i\in[n]$.
If $x^{1}_i(\tau)=0$, then from \rep{updates} and Assumption \ref{para}, $\dot x^{1}_i(\tau)\ge 0$.
The same holds for $x^{2}_i(\tau)$ and $x^{1}_i(\tau)+x^{2}_i(\tau)$.
If $x^{1}_i(\tau)=1$, then from \rep{updates} and Assumption \ref{para}, $\dot x^{1}_i(\tau)\le 0$.
The same holds for $x^{2}_i(\tau)$ and $x^{1}_i(\tau)+x^{2}_i(\tau)$.
It follows that $x^{1}_i(t),x^{2}_i(t),x^{1}_i(t)+x^{2}_i(t)$ will be in $[0,1]$ for all times $t\ge \tau$.
Since the above arguments hold for all $i\in[n]$,  $x^{1}_i(t),x^{2}_i(t),x^{1}_i(t)+x^{2}_i(t)$ will be in $[0,1]$ for
all $i\in[n]$ and $t\ge \tau$.
Since by Assumption \ref{x0} $x^{1}_i(0),x^{2}_i(0),x^{1}_i(0)+x^{2}_i(0)\in[0,1]$ for all $i\in[n]$,
it follows that $x^{1}_i(t),x^{2}_i(t),x^{1}_i(t)+x^{2}_i(t)\in[0,1]$ for all $i\in[n]$ and $t\ge 0$.
\hfill
$\qed$

{\em Proof of Proposition \ref{single0} (case $s(-D+B)\le 0$):}
We first consider the case when  $s(-D+B)<0$.
Since $(-D+B)$ is an irreducible Metzler matrix, by Lemma \ref{less},
there exists a positive diagonal matrix $P$ such that
$(-D+B)'P+P(-D+B)$ is negative definite.
Consider the Lyapunov function candidate  $V(x(t))=x(t)'Px(t)$.
From \rep{single}, when $x(t)\neq \0$,
\begin{eqnarray*}
\dot V(x(t)) &=& 2x(t)'P\left(-D+B-X(t)B\right)x(t) \\
&< & -2 x(t)'PX(t)Bx(t) \leq  0.
\end{eqnarray*}
Thus, in this case, $\dot V(x(t))<0$ if $x(t)\neq \0$.
By Lemma \ref{box} and Corollary \ref{lya},
$x=\0$ is asymptotically stable with  domain of attraction $[0,1]^n$.

Next we consider the case when  $s(-D+B)=0$.
Since $(-D+B)$ is an irreducible Metzler matrix, by Lemma \ref{equal},
there exists a positive diagonal matrix $P$ such that
$(-D+B)'P+P(-D+B)$ is negative semi-definite.
Consider the Lyapunov function candidate $V(x(t))=x(t)'Px(t)$.
From \rep{single}, we have
\begin{eqnarray*}
\dot V(x(t)) &=& 2x(t)'P\left(-D+B-X(t)B\right)x(t) \\
&= & x(t)'\left((-D+B)'P+P(-D+B)\right)x(t) \\
&& \ \ \ \  \ \ \ \ \  \ - 2x(t)'PX(t)Bx(t) \leq 0.
\end{eqnarray*}
We claim that $\dot V(x(t))<0$ if $x(t)\neq \0$.
To establish this claim, we first consider the case when $x(t)\gg \0$.
Since $B$ is nonnegative and irreducible, $Bx(t)\gg \0$.
Since $P$ is a positive diagonal matrix, it follows that $x(t)'PX(t)Bx(t)>0$,
so $\dot V(x(t))<0$.
Next we consider the case when $x(t)>\0$ and $x(t)$ has at least one zero entry.
Since $(-D+B)$ is an irreducible Metzler matrix and $P$ is a positive diagonal matrix,
$(-D+B)'P+P(-D+B)$ is a symmetric irreducible Metzler matrix.
Since $(-D+B)'P+P(-D+B)$ is negative semi-definite,
it follows that $s((-D+B)'P+P(-D+B))=0$. By Lemma \ref{metzler0}, $0$ is a simple eigenvalue of
$(-D+B)'P+P(-D+B)$ and it has a unique (up to scalar multiple) strictly positive eigenvector
corresponding to the eigenvalue $0$. That is, there exists $\tilde{x}\gg 0$ such that $A\tilde{x} = 0 \tilde{x}$. Thus, $x(t)'\left((-D+B)'P+P(-D+B)\right)x(t)<0$ when
$x(t)>\0$ and $x(t)$ has at least one zero entry. 
Therefore, $\dot V(x(t))<0$ if $x(t)\neq \0$.
By Lemma \ref{box} and Corollary \ref{lya},
$x=\0$ is asymptotically stable with  domain of attraction $[0,1]^n$.
\hfill
$\qed$

To proceed, we need the following lemma.

\begin{lemma}
Suppose that $\delta_i\ge0$ for all $i\in[n]$ and that matrix $B$ is nonnegative and irreducible.
If $x^*$ is a nonzero equilibrium of system~\rep{single}, then $x^*\gg \0$.
\label{0}
\end{lemma}

{\em Proof of Lemma~\ref{0}:}
Suppose that $x^*$ is a nonzero equilibrium of \rep{single}.
By Lemma \ref{box}, it must be true that
$x^*\geq \0$.
To prove $x^*\gg \0$, suppose that, to the contrary, $x^*$ has at least one zero entry.
Without loss of generality, set $x^*_1=0$. Since $x^*$ is an equilibrium of \rep{single}, from \rep{update},
$$-\delta_1x^*_1 + (1-x^*_1)\sum_{j=1}^n \beta_{1j}x^*_j = \sum_{j=1}^n \beta_{1j}x^*_j = 0.$$
It then follows that for any $j\in[n]$ such that $\beta_{1j}>0$, $x^*_j=0$.
By repeating this argument, since $B$ is irreducible,
we have $x^*_i=0$ for all $i\in[n]$. This contradicts the assumption that $x^*\geq 0$.
Thus, $x^*\gg \0$.
\hfill
$\qed$

{\em Proof of Proposition \ref{single0} (case $s(-D+B)> 0$):}
It will suffice to show that if $s(-D+B)> 0$, then there exists a unique strictly positive equilibrium.
We first show that there exists an $x^*\gg \0$ which is an equilibrium of \rep{single}.

Let $c>0$ be any positive constant such that
\eq{s(-D+B) -c > 0. \label{c}}
Such a constant $c$ always exists since $s(-D+B)> 0$.
Set
$\bar D = D+cI$.
Since $\delta_i\ge0$ for all $i\in[n]$, 
$D$ is a nonnegative diagonal matrix.
Thus, $\bar D$ is nonsingular and $\bar D^{-1}$ is also a positive diagonal matrix.

Consider the above equation and define a continuous map $f:(0,1]^n\rightarrow[0,1]^n$ given by
$$f(x) = \left(I-c\bar D^{-1}+{\rm diag}(\bar D^{-1}Bx)\right)^{-1}\bar D^{-1}Bx.$$
Since the domain of $f$ is $(0,1]^n$, $x$ as the argument of $f$ satisfies $x\gg \0$. Therefore ${\rm diag}(\bar D^{-1}Bx)$ is a positive diagonal matrix, and we have that $(I-c\bar D^{-1}+{\rm diag}(\bar D^{-1}Bx))$ is invertible. Therefore $f$ is well-defined. 
Note that the $i$th entry of $f(x)$, denoted by $f_i(x)$, is given by
$$f_i(x) = \frac{\left(\bar D^{-1}Bx\right)_i}{1-\frac{c}{c+\delta_i}+\left(\bar D^{-1}Bx\right)_i}.$$ 
Since $\bar D^{-1}$ and $B$ are both nonnegative,
for any $y \geq z$ in $(0,1]^n$,  $f_i(y)\geq f_i(z)$,
so $f(y)\geq f(z)$.

Since $\bar D^{-1}B$ is an irreducible nonnegative matrix, by Lemma \ref{metzler0},
there exists $v\gg \0$ such that
\eq{\bar D^{-1}Bv = rv, \label{r}}
where
$r=\rho(\bar D^{-1}B)$.
Since $s(-\bar D +B) = s(-D+B)-c $, from \rep{c}, $s(-\bar D +B)>0$.
By Proposition \ref{ff}, it follows that
$r>1$. Then, we can always find an $\varepsilon>0$
such that for each $i\in[n]$,
\eq{\varepsilon v_i \le \frac{r-1}{r}.\label{small}}
From this, it follows that
$1\le \frac{r}{1+\varepsilon rv_i}$,
and thus,
$\varepsilon v_i \le \frac{\varepsilon rv_i}{1+\varepsilon rv_i}$.
From \rep{r}, we have 
$$\varepsilon v_i \le \frac{\left(\bar D^{-1}B\varepsilon v\right)_i}{1+\left(\bar D^{-1}B\varepsilon v\right)_i}
\le \frac{\left(\bar D^{-1}B\varepsilon v\right)_i}{1-\frac{c}{c+\delta_i}+\left(\bar D^{-1}B\varepsilon v\right)_i},$$
which implies that $\varepsilon v \le f(\varepsilon v)$.
It follows from \rep{small} that $\varepsilon v \ll \1$.
Since we have already shown that for any $y \geq z$ in $(0,1]^n$, 
$f(y)\geq f(z)$,
$f$ maps the compact, convex set $\scr{C}=\{x\ | \ \varepsilon v \le x \le \mathbf{1} \}$
to itself. By Brouwer's fixed-point theorem, $f$ has a fixed point in $\scr{C}$,
which must be strictly positive.
Let $x^*\gg \0$ denote this fixed point. Then $f(x^*) = x^*$, i.e., 
$$x^*=\left(I-c\bar D^{-1}+{\rm diag}(\bar D^{-1}Bx^*)\right)^{-1}\bar D^{-1}Bx^*.$$
Therefore, we have 
\begin{align*}
\bar D^{-1}Bx^* &= \left(I-c\bar D^{-1}+{\rm diag}(\bar D^{-1}Bx^*)\right)x^*\\
&= x^*+{\rm diag}(\bar D^{-1}Bx^*)x^*-c\bar D^{-1}x^*\\
&= x^*+X^*\bar D^{-1}Bx^*-c\bar D^{-1}x^* \\
&= x^*+\bar D^{-1}X^*Bx^*-c\bar D^{-1}x^*.
\end{align*}
Therefore, 
$Bx^*=\bar Dx^*+X^*Bx^* - cx^* = Dx^*+X^*Bx^*$,
by definition of $\bar D$. 
Thus, we have 
$((-D+B)-X^*B)x^*=\0$,
and therefore $x^*\gg \0$ is an equilibrium of \rep{single}. 

It remains to show that the strictly positive equilibrium  
is unique.
Suppose that $x$ and $y$ are both nonzero equilibria of \rep{single}, and let $\varepsilon$ from \rep{small} be sufficiently small such that $x,y \in \scr{C}$. From Lemma \ref{0},
it follows that $x,y\gg \0$. Set
$$\kappa = \max_{i\in [n]} \frac{x_i}{y_i}.$$
Then, $x\le \kappa y$, and there exists $j\in[n]$ for which $x_j = \kappa y_j$.
We claim that $\kappa \le 1$. To establish this claim, suppose that,
to the contrary, $\kappa >1$.
Since $x$ is a fixed point of $f$ and for any $u \geq v$ in $(0,1]^n$, 
$f_j(u)\geq f_j(v)$ for all $j\in[n]$,
it follows that
\begin{eqnarray*}
x_j &=& \frac{\left(\bar D^{-1}Bx\right)_j}{1-\frac{c}{c+\delta_j}+\left(\bar D^{-1}Bx\right)_j} \\
&\le&  \frac{\left(\bar D^{-1}B\kappa y\right)_j}{1-\frac{c}{c+\delta_j}+\left(\bar D^{-1}B\kappa y\right)_j} \\
&=&
\frac{\kappa\left(\bar D^{-1}B y\right)_j}{1-\frac{c}{c+\delta_j}+\kappa\left(\bar D^{-1}B y\right)_j}.
\end{eqnarray*}
From the assumption that $\kappa >1$, we have
$$\frac{\kappa\left(\bar D^{-1}B y\right)_j}{1-\frac{c}{c+\delta_j}+\kappa\left(\bar D^{-1}B y\right)_j}<
\frac{\kappa\left(\bar D^{-1}B y\right)_j}{1-\frac{c}{c+\delta_j}+\left(\bar D^{-1}B y\right)_j}
.$$
Since $y$ is a fixed point of $f$,
$$
\frac{\left(\bar D^{-1}B y\right)_j}{1-\frac{c}{c+\delta_j}+\left(\bar D^{-1}B y\right)_j}
= y_j.$$
Then, it follows that
$$x_j <
\frac{\kappa\left(\bar D^{-1}B y\right)_j}{1-\frac{c}{c+\delta_j}+\left(\bar D^{-1}B y\right)_j} = \kappa y_j =x_j,$$
which is a contradiction. Therefore, $\kappa \le 1$,
which implies that $x\le y$. Using the same arguments and by exchanging the roles of $x$ and $y$,
it also can be shown that $y\le x$. Thus, $x=y$, which establishes
the uniqueness of the strictly positive equilibrium.
This completes the proof.
\hfill
$\qed$

{\em Proof of Lemma~\ref{attract}:}
If $x(0)\gg \0$, then the lemma is true with $\tau=0$.
Suppose that $x(0)> \0$ with $x_i(0)=0$ for at least one $i\in[n]$.
Let $\scr{F}(t)$ be the set of all $i\in[n]$ such that
$x_i(t)=0$. 
In other words, $x_i(t)=0$ for all $i\in\scr{F}(t)$
and $x_i(t)>0$ for all $i\in[n]\setminus\scr{F}(t)$.
Clearly, the set $\scr{F}(0)$ is nonempty. 
Moreover, since the matrix $B$ is irreducible, 
there exists $j\in\scr{F}(0)$ such that
$x_j(0)=0$, $x_k(0)>0$, and $\beta_{jk}>0$.
From \rep{update}, it follows that $\dot x_j(0)>0$.
Thus, there must exist $\tau_1>0$ such that $x_j(\tau_1)>0$
and $x_i(\tau_1)>0$ for all $i\in[n]\setminus\scr{F}(0)$.
This implies that $\scr{F}(\tau_1)$ is a proper subset of $\scr{F}(0)$. 
Note that $\scr{F}(0)$ is a finite set. By repeating the above arguments, we conclude that 
there exists $\tau>0$ such that $\scr{F}(\tau)$ is the empty set,
which implies that $x_i(\tau)>0$ for all $i\in[n]$.
\hfill
$\qed$

{\em Proof of Proposition \ref{2local}:}
Let $y_i(t) = x_i(t)- x^*_i$ for all $i\in[n]$, i.e.,  $y(t)=x(t)-x^*$. 
Let $Y(t)={\rm diag}(y(t))$ and $X^*={\rm diag}(x^*)$.
Note that
$ \left(-D+B-X^*B\right)x^* =\0$.
Then, 
\begin{eqnarray*}
\dot y(t) &=& \left( -D+B-\left(Y(t)+X^*\right)B \right)\left(y(t)+x^*\right)\\
&=& \left( -D+\left(I-X^*\right)B -Y(t)B \right)y(t) - Y(t)Bx^* \\
&=& \left( -D+\left(I-X^*\right)B \right)y(t) - Y(t)Bx(t) \\
&=& \left( -D+\left(I-X^*\right)B - {\rm diag}(Bx(t))\right)y(t).
\end{eqnarray*}
Thus, for all $i\in[n]$,
$$\dot y_i(t) = -\delta_i y_i(t) + (1-x^*_i)\sum_{j=1}^n \beta_{ij}y_j(t) - \left(\sum_{j=1}^n \beta_{ij}x_j(t)\right)y_i(t).$$

Consider the Lyapunov function candidate 
$$V(y(t))=\max_{k\in[n]}\frac{|y_k(t)|}{x^*_k},$$
which is well-defined since $x^*_k >0$ for all $k\in[n]$  by Proposition \ref{single0}. Note that $V(y(t))\ge 0$, with equality if and only if $y(t)=\0$ (or equivalently, $x(t)=x^*$).
Moreover, for all $i\in[n]$ and $t\geq 0$, there holds
$|y_i(t)| \le V(y(t)) x^*_i$.
For any time $t\geq 0$, without loss of generality, let $m\in[n]$ such that 
$$\frac{|y_m(t)|}{x^*_m}=V(y(t))=\max_{k\in[n]}\frac{|y_k(t)|}{x^*_k}.$$
Then, when $|y_m(t)|>0$ (or equivalently, $V(y(t))>0$),

$ \displaystyle\dot V(y(t)) = \frac{1}{x^*_m}\cdot \frac{d|y_m(t)|}{dt}$
\begin{eqnarray}
&= & \frac{1}{x^*_m} {\rm sgn}(y_m(t))\dot y_m(t) \label{kkk}\\
&= & \frac{1}{x^*_m} {\rm sgn}(y_m(t))\Bigg( -\delta_m y_m(t) + (1-x^*_m)\sum_{j=1}^n \beta_{mj}y_j(t) \nonumber \\
&&\ \ \ \ \ \ \ \ \ \ \ \ \ \ \ \ \ \ \ \ \ \ \ \ - \Bigg(\sum_{j=1}^n \beta_{mj}x_j(t)\Bigg)y_m(t)\Bigg) \nonumber\\
&=& \frac{1}{x^*_m} \Bigg( -\delta_m |y_m(t)| + (1-x^*_m)\sum_{j=1}^n \beta_{mj}y_j(t){\rm sgn}(y_m(t)) \nonumber \\
&&\ \ \ \ \ \ \ \ \ \ \ \ \ \ \ \ \ \ \ \ \ \ \ \ - \Bigg(\sum_{j=1}^n \beta_{mj}x_j(t)\Bigg)|y_m(t)|\Bigg) \nonumber\\
&\le & \frac{1}{x^*_m} \Bigg( -\delta_m |y_m(t)| + (1-x^*_m)\sum_{j=1}^n \beta_{mj}|y_j(t)|\Bigg) \nonumber \\
&&\ \ \ \ \ \ \ \ \ \ \ \ \ \ \ \ \ \ \ \ \ \ \ \ - \frac{1}{x^*_m} \Bigg(\sum_{j=1}^n \beta_{mj}x_j(t)\Bigg)|y_m(t)| \nonumber \\
&\le & \frac{V(y(t))}{x^*_m} \Bigg( -\delta_m x^*_m + (1-x^*_m)\sum_{j=1}^n \beta_{mj}x^*_j\Bigg) \nonumber \\
&&\ \ \ \ \ \ \ \ \ \ \ \ \ \ \ \ \ \ \ \ \ \ \ \ - \frac{1}{x^*_m} \Bigg(\sum_{j=1}^n \beta_{mj}x_j(t)\Bigg)|y_m(t)| \nonumber \\
&=& - \frac{1}{x^*_m} \Bigg(\sum_{j=1}^n \beta_{mj}x_j(t)\Bigg)|y_m(t)| \label{jjj}\\
&\le & 0. \nonumber
\end{eqnarray}
\normalsize
From \rep{kkk} and the definition of $y_i(t)$, it is straightforward to verify that
$\dot V(y(t)) = 0$ in the case when $x(t)=x^*$ (i.e., $y(t)=\0$).
Next we consider the case when $x(t)\gg \0$.
From \rep{jjj}, since the matrix $B$ is irreducible  and
$x^* \gg \0$ by Proposition \ref{single0}, it can be seen that $\dot V(y(t)) < 0$
if $x(t)\neq x^*$.

Recall that by Lemma \ref{attract}, as long as $x(0)> \0$, there always exists
a finite time $t_0$ at which $x(t_0)\gg \0$.
From Lemma \ref{box} (with $x^2 = \0$, this reduces to the single virus model case) and Corollary \ref{lya}, to prove the proposition,
it remains to show that the system \rep{single} is positively invariant on a subset of $(0,1]^n$.
Without loss of generality,
suppose that $x(0)\gg \0$.
Let $\varepsilon$ be a nonnegative real number such that
$$\varepsilon = \max_{k\in[n]}\frac{|x_k(0)-x^*_k|}{x^*_k}.$$
Consider the set
$\scr{B}= \left\{ x \ | \ V(x(t)) \le \varepsilon\right\} \cap \ (0,1]^n$.
Since $\dot V(x(t)) \le 0$, the system~\rep{single} is positively invariant on $\{ x \ | \ V(x(t)) \le \varepsilon\}$. 
Suppose that there exists some finite time $t_1$ at which $x(t_1)>\0$ and $x(t_1)\notin (0,1]^n$. From Lemma~\ref{attract}, there must exist another finite time $t_2>t_1$ such that $x(t_2)\gg \0$. Since we have shown that $\dot V(x(t)) < 0$ if $x(t)\gg \0$ and $x(t)\neq x^*$, it follows that $x(t)\in (0,1]^n$ for all $t\ge t_2$. 
Note that
$$V(x(t))=\max_{k\in[n]}\frac{|x_k(t)-x^*_k|}{x^*_k}.$$
It then follows that for any $z\in (0,1]^n$ which is close to $\0$, there holds $V(z)< V(\0)$. Since $\dot V(x(t)) < 0$ if $x(t)\gg \0$ and $x(t)\neq x^*$,
there cannot exist a trajectory from $x(0)$ to $\0$.
Thus, the system~\rep{single} is positively invariant on $(0,1]^n$ for $t\ge t_2$,
which implies that the system~\rep{single} is positively invariant on $\scr B$ for $t\ge t_2$.
Since the above arguments hold for any $x(0)\gg \0$,
such a set $\scr{B}$ always exists, and thus after some finite time, 
the system \rep{single} is positively invariant on $\scr{B}\subset (0,1]^n$.
\hfill
$\qed$

{\em Proof of Lemma~\ref{nonzero}:}
Let $(\tilde x^{1},\tilde x^{2})$ be an equilibrium of \rep{sys} and suppose, to the contrary, that there exists some $i\in[n]$ such that
$\tilde x^{1}_i+\tilde x^{2}_i = 1$. 
Let $(\tilde x^{1}(0),\tilde x^{2}(0)) = (\tilde x^{1},\tilde x^{2})$.
Then, from \rep{updates},  $\dot{\tilde x}^{1}_i(0)<0$ and/or $\dot{\tilde x}^{2}_i(0) < 0$, which contradicts
the hypothesis that $(\tilde x^{1}_i,\tilde x^{2}_i)$ is an equilibrium. Therefore, $\tilde x^{1}(\tau)+\tilde x^{2}(\tau)\ll \1$.
\hfill
$\qed$


%

\bibliographystyle{IEEEtran}
\bibliography{bib}


\begin{IEEEbiography} [{\includegraphics[width=1in,height=1.25in,clip,keepaspectratio]{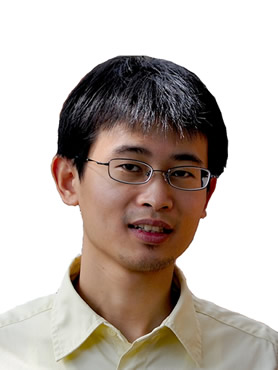}}]
{Ji Liu} received the B.S. degree in information engineering from Shanghai Jiao Tong University, Shanghai, China, in 2006, 
and the Ph.D. degree in electrical engineering from Yale University, New Haven, CT, USA, in 2013. He is currently an 
Assistant Professor in the Department of Electrical and Computer Engineering at Stony Brook University, Stony Brook, NY, 
USA. Prior to joining Stony Brook University, he was a Postdoctoral Research Associate at the Coordinated Science 
Laboratory, University of Illinois at Urbana-Champaign, Urbana, IL, USA, and the School of Electrical, Computer and Energy 
Engineering, Arizona State University, Tempe, AZ, USA. His current research interests include distributed control and
computation, distributed optimization and learning, multi-agent systems, social networks, epidemic networks, and 
cyber-physical systems.
\end{IEEEbiography}


\begin{IEEEbiography} [{\includegraphics[width=1in,height=1.25in,clip,keepaspectratio]{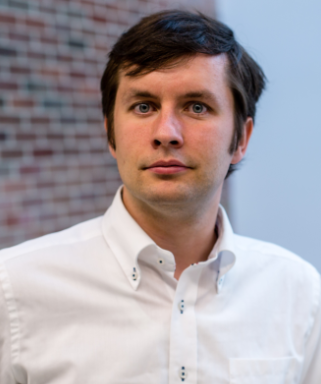}}]
{Philip E. Par\'{e}} received his B.S. in Mathematics with University Honors and his M.S. in Computer Science from Brigham Young University, Provo, UT, in 2012 and 2014, respectively, and his Ph.D. in Electrical and Computer Engineering from the University of Illinois at Urbana-Champaign, Urbana, IL in 2018.  He was the recipient of the 2017-2018 Robert T. Chien Memorial Award for excellence in research and named a 2017-2018 College of Engineering Mavis Future Faculty Fellow. 
His research interests include the modeling and control of dynamic networked systems, model reduction techniques, and time--varying systems.
\end{IEEEbiography}


\begin{IEEEbiography} [{\includegraphics[width=1in,height=1.25in,clip,keepaspectratio]{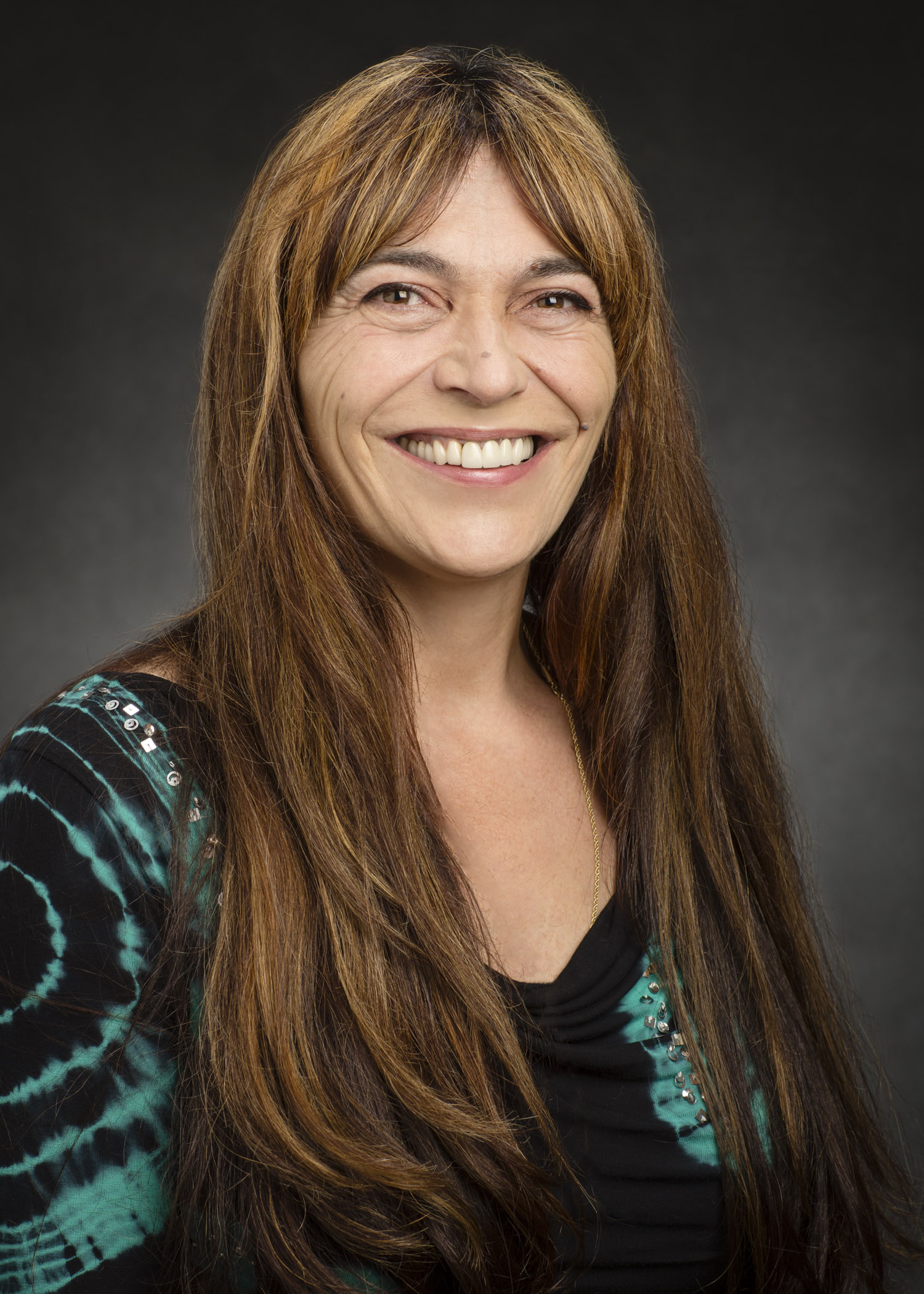}}]
{Angelia Nedi\'{c}} has a Ph.D. from Moscow State University, Moscow, Russia, 
in Computational Mathematics and Mathematical Physics (1994), and a Ph.D. 
from Massachusetts Institute of Technology, Cambridge, USA in Electrical 
and Computer Science Engineering (2002). She has worked as a senior engineer 
in BAE Systems North America, Advanced Information Technology Division at 
Burlington, MA. Currently, she is a faculty member of the school of Electrical, 
Computer and Energy Engineering at Arizona State University at Tempe. 
Prior to joining Arizona State University, she was a Willard Scholar 
faculty member at the University of Illinois at Urbana-Champaign. She has 
been a recipient of NSF CAREER Award 2007 in Operations Research for her work 
in distributed multi-agent optimization. She is a recipient (jointly with her co-authors) 
of the Best Paper Award at the Winter Simulation Conference 2013 and the Best Paper Award 
at the International Symposium on Modeling and Optimization in Mobile, Ad Hoc 
and Wireless Networks (WiOpt) 2015. Her general research interest is in large scale 
complex systems dynamics and optimization.
\end{IEEEbiography}


\begin{IEEEbiography} [{\includegraphics[width=1in,height=1.25in,clip,keepaspectratio]{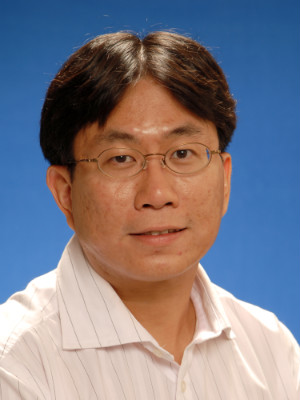}}]
{Choon Yik Tang} (S'97-M'04) received the B.S. and M.S. degrees in mechanical engineering from Oklahoma State University, Stillwater, in 1996 and 1997, respectively, and the Ph.D. degree in electrical engineering from the University of Michigan, Ann Arbor, in 2003. He was a Postdoctoral Research Fellow at the University of Michigan from 2003 to 2004, a Research Scientist at Honeywell Labs, Minneapolis, from 2004 to 2006, and a Visiting Scholar at the University of Illinois at Urbana-Champaign in 2014. Since 2006, he has been with the School of Electrical and Computer Engineering, University of Oklahoma, Norman, where he is currently an Associate Professor. His research interests include systems and control theory, distributed algorithms for computation and optimization over networks, and control and operation of wind farms.
\end{IEEEbiography}


\begin{IEEEbiography} [{\includegraphics[width=1in,height=1.25in,clip,keepaspectratio]{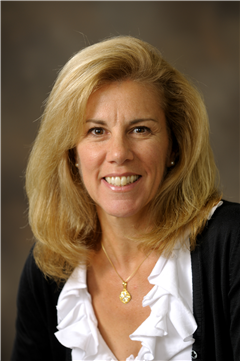}}]
{Carolyn L. Beck} is an associate professor at the University of Illinois, Urbana-Champaign (UIUC) in the Industrial and Enterprise Systems Engineering (IESE) Department. She completed her Ph.D. at Caltech, her M.S. at Carnegie Mellon, and her B.S. at Cal Poly, all in Electrical Engineering. Prior to completing her Ph.D., she worked at Hewlett-Packard in Silicon Valley for four years, designing digital hardware and software for measurement instruments.  She has held visiting faculty positions at the Royal Institute of Technology (KTH) in Stockholm, Stanford University in Palo Alto and Lund University in Lund, Sweden.   She has received national research awards and local teaching awards.  Prof. Beck's research interests range from network inference problems to control of anesthetic pharmacodynamics. Her main research interests are: model reduction and approximation for the purpose of feedback control design; mathematical systems theory; clustering and aggregation methods. 
\end{IEEEbiography}


\begin{IEEEbiography} [{\includegraphics[width=1in,height=1.25in,clip,keepaspectratio]{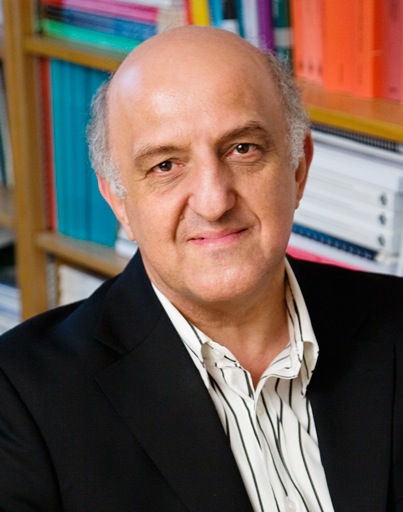}}]
{\bf Tamer Ba\c{s}ar} (S'71-M'73-SM'79-F'83-LF'13) is with the University of Illinois at Urbana-Champaign (UIUC), where he holds the academic positions of
Swanlund Endowed Chair;
Center for Advanced Study Professor of  Electrical and Computer Engineering;
Research Professor at the Coordinated Science
Laboratory; and Research Professor  at the Information Trust Institute.
He is also the Director of the Center for Advanced Study.
He received B.S.E.E. from Robert College, Istanbul,
and M.S., M.Phil, and Ph.D. from Yale University. He is a member of the US National Academy
of Engineering,  and Fellow of IEEE, IFAC and SIAM, and has served as president of IEEE CSS,
ISDG, and AACC. He has received several awards and recognitions over the years, including the
highest awards of IEEE CSS, IFAC, AACC, and ISDG; the IEEE Control Systems Award; and a number of international honorary
doctorates and professorships. He has over 900 publications in systems, control, communications,
and dynamic games, including books on non-cooperative dynamic game theory, robust control,
network security, wireless and communication networks, and stochastic networked control. He was
the Editor-in-Chief of Automatica between 2004 and 2014, and is currently  editor of several book series. His current research interests
include stochastic teams, games, and networks; security; and cyber-physical systems.
\end{IEEEbiography}

\end{document}